\input ./style/arxiv-general.cfg
\documentclass[aap,MSNbibl,seceqn,dvips]{arximspdf}
\makeatletter
   \@ifpackageloaded{graphicx}{}{\usepackage{graphicx}}
\makeatother
\usepackage{mathbh}

\doi{10.1214/15-AAP1141}
\volume{26}
\issue{4}
\pubyear{2016}
\firstpage{2106}
\lastpage{2140}
\docsubty{FLA}

\makeatletter

\def\vafrac#1#2{(#1)/(#2)}

\newcommand{\mrr}{\mathrm{d}}
\newcommand{\rrvert}{\vert}
\newcommand{\rrVert}{\Vert}
\newcommand{\llvert}{\vert}
\newcommand{\llVert}{\Vert}
\renewcommand{\mid}{|}
\newcommand{\eqref}[1]{(\ref{#1})}
\newcommand{\R}{\mathbb{R}}
\newtheorem{Theorem}{Theorem}[section]
\newtheorem{Proposition}[Theorem]{Proposition}
\newtheorem{Corollary}[Theorem]{Corollary}
\newtheorem{Lemma}[Theorem]{Lemma}
\newproclaim{Remark}[Theorem]{Remark}
\newproclaim{Example}[Theorem]{Example}
\def\ind{\mathbh{1}}
\makeatother

\begin{document}
\begin{frontmatter}

\title{Gaussian approximation of nonlinear Hawkes~processes}
\runtitle{Approximation of Hawkes processes}

\begin{aug}
\author[A]{\fnms{Giovanni Luca}~\snm{Torrisi}\corref{}\ead
[label=e2]{torrisi@iac.rm.cnr.it}}
\runauthor{G. L. Torrisi}
\affiliation{CNR---Istituto per le Applicazioni del Calcolo ``Mauro Picone''}
\address[A]{CNR---Istituto per le Applicazioni\\
\quad del Calcolo ``Mauro Picone''\\
Via dei Taurini 19\\
I-00185 Rome\\
Italy\\
\printead{e2}}
\end{aug}
%

%
\received{\smonth{10} \syear{2014}}
%
\revised{\smonth{6} \syear{2015}}

%
\begin{abstract}
We give a general Gaussian bound for the first chaos (or innovation) of
point processes with stochastic intensity constructed by embedding in a bivariate Poisson
process. We apply the general result to nonlinear Hawkes processes, providing quantitative
central limit theorems.
\end{abstract}

%
\end{frontmatter}

\section{Introduction}\label{sec:int}

In the seminal papers \cite{nourdinpeccati} and \cite{peccati},
Stein's method and Malliavin's calculus have been combined to derive
explicit bounds in the
Gaussian approximation of random variables on the Wiener and Poisson
spaces. Further developments on the Poisson space include, for example,
\cite{last0,peccatizheng,privaulttorrisi,reitzner}.
In particular, in \cite{privaulttorrisi} the authors derive new
Gaussian bounds
for functionals of the one-dimensional homogeneous Poisson process by
using the Clark--Ocone representation formula; see, for example, \cite
{privault}.
In contrast with covariance identities based
on the inverse of the Ornstein--Uhlenbeck operator the Clark--Ocone
representation formula only requires the computation of a gradient
and a conditional expectation. For this reason, the Clark--Ocone
representation formula is a valuable tool even for the probability approximation
of random variables on spaces different from the Wiener and Poisson. We
refer the reader to \cite{privaulttorrisialea} for the
use of the Clark--Ocone representation formula for the Gaussian and
Poisson approximation of random variables on the Bernoulli space and to
\cite{nguyenprivaulttorrisi}
for the use of the Clark--Ocone representation formula for the Gaussian
approximation of solutions of some stochastic equations.

The contributions of this paper are the following. We provide a
Gaussian bound for the first chaos of
a large class of point processes with stochastic intensity; see Theorem
\ref{thm:normalskor}. Particularly, we consider
point processes on the line constructed by embedding in a bivariate
Poisson process and provide a Gaussian approximation for the first
chaos (or innovation)
combining Stein's method and Malliavin's calculus via a Clark--Ocone
type representation formula.

To the best of our knowledge, this is the first paper which provides
Gaussian bounds
for the innovation of a point process with stochastic intensity by the
Malliavin--Stein method.

We apply our general result to nonlinear Hawkes processes,
deriving an explicit Gaussian bound for the innovation; see Theorem
\ref{thm:nonlinGauss}. In the special case of self-exciting processes
(or linear Hawkes processes),
relying on the knowledge of the intensity of the process and the
spectral theory of point processes, we are able to provide alternative
Gaussian bounds for the innovation
which, in some cases, improve those one obtained by directly applying
Theorem \ref{thm:nonlinGauss}; see Theorems \ref{thm:linGauss1}, \ref
{thm:linGaussappvecchio} and Proposition
\ref{prop:compare}. We exploit such Gaussian bounds to provide new
quantitative central limit theorems in the Wasserstein distance for the
first chaos of Hawkes processes;
see Corollaries \ref{cor:QCLTnonlinear} and \ref{cor:QCLT}. The
quantitative nature of these Gaussian approximations allows, for example,
to construct in a standard way confidence intervals for the
corresponding innovations, we outlined this simple application in
Example \ref{ex:1}.

From the point of view of applications, the extension of our results to
multivariate point processes with stochastic intensity and random marks
is certainly of interest;
see, for example, \cite{bacry1,massoulie} and \cite
{torrisi}. This topic is presently under investigation by the author,
as well as the topic concerning the Poisson approximation, via the
Malliavin--Stein method,
of first-order stochastic integrals with respect to point processes
with stochastic intensity (note that for this latter argument
some results are already known, see \cite{BBAOP} and \cite{BBSPA}).

In the last years, there has been a renewed interest on Hawkes
processes, mainly due to their
mathematical tractability and versatility in modeling contexts.
Self-exciting processes were introduced in \cite{hawkes} and \cite{hawkes1},
while the wider class of nonlinear Hawkes processes was introduced in
\cite{bremaudmassoulie1}. Various mathematical aspects
of these processes (and their generalizations), such as stability, rate
of convergence to equilibrium, perfect and approximate simulation,
large deviations and limit theorems,
are studied in \cite
{bordenave,bremaudmassoulie1,bremaudmassoulie,bnt,hawkes,hawkes1,jaisson,massoulie,moller1,moller2,moller3,torrisi,zhu1,zhu2}.
Linear Hawkes processes are Poisson cluster processes with a simple
self-exciting structure which makes them very appealing to account for
situations
where the occurrence of future events directly depends on the past
history. Nonlinear Hawkes processes allow to account for inhibitory effects.
For these reasons, Hawkes processes naturally and simply capture a
causal structure
of discrete events dynamics associated with endogenous triggering,
contagion and self-activation phenomena. Typical fields where this kind
of dynamics arise are
seismology (occurrence of earthquakes), neuroscience (occurrence of
neuron's spikes), genome analysis (occurrence of events along a DNA
sequence), insurance (occurrence of claims) and finance (occurrence of
market order arrivals);
see, for example, \cite
{bacry1,bacry2,kagan,ogata,reino1,reino2,stabile,zhu0} for applications
of Hawkes processes
in these contexts.

The paper is organized as follows. In Section~\ref{sec:2}, we give
some preliminaries on point processes including the notion of
stochastic intensity,
the Poisson embedding construction and a Clark--Ocone type
representation formula. In Section~\ref{sec:3},
we prove a general upper bound for the Wasserstein distance between the
first chaos of a point process with stochastic intensity
(constructed by embedding on a bivariate Poisson process) and a
standard normal random variable. In Section~\ref{sec:4}, we apply the result
in Section~\ref{sec:3} to nonlinear Hawkes processes. Particularly, in
Section~\ref{subsec:4/1} we provide an explicit Gaussian bound
for the first chaos (and a suitable approximated version of it) of a
stationary nonlinear Hawkes process. The corresponding
quantitative central limit theorem is derived in Section~\ref
{subsec:4/2}. The special case of self-exciting processes
is treated in Section~\ref{sec:linear}.

\section{Preliminaries on point processes}\label{sec:2}

In this section, we give some preliminaries on point processes, and
refer the reader to the books \cite{bremaud,daley,daley2} for
more insight into this subject.

Let $\{T_n\}_{n\in\mathbb{Z}}$ be a sequence of random times defined
on a probability space $(\Omega,\mathcal{A},P)$. Given a Borel set
$A\in\mathcal{B}(\R)$, we define
\[
N(A):=\sum_{n\in\mathbb{Z}}\ind_{A}(T_n)
\]
and we call $N:=\{N(A)\}_{A\in\mathcal{B}(\R)}$ the point process
with times $\{T_n\}_{n\in\mathbb{Z}}$. We suppose that $N$ has the
following properties:
\begin{eqnarray*}
T_n\in\overline{\R}&:=&\R\cup\{\pm\infty\};\qquad \llvert
T_n\rrvert <\infty\quad\Longrightarrow\quad T_n<T_{n+1};
\qquad 
T_0\leq0<T_1;
\\
N(A)&<&\infty,\qquad\mbox{for all bounded $A$.} 
\end{eqnarray*}
%
These conditions guarantee that $N$ is simple, that is, $N(\{a\})\leq1$
for any $a\in\R$, and locally finite.

Given a sequence $\{Z_n\}_{n\in\mathbb{Z}}$ of random variables on
$\Omega$ with values in some measurable space $(E,\mathcal{E})$, we define
\[
\overline{N}(A):=\sum_{n\in\mathbb{Z}}\ind_{A}(T_n,Z_n),
\qquad A\in\mathcal{B}(\R)\otimes\mathcal{E}
\]
and
\[
\int_A\psi(t,z)\overline{N}(\mrr t\times \mrr z):=\sum_{n\in\mathbb{Z}}\psi(T_n,Z_n)
\ind_{A}(T_n,Z_n)
\]
for a measurable function $\psi:\R\times E\to\R$ for which the
infinite sum is well defined.

\subsection{Point processes with stochastic intensity}

Let $\mathcal{F}:=\{\mathcal{F}_t\}_{t\in\R}\subset\mathcal{A}$
be a filtration such that
$\mathcal{F}_t\supseteq\mathcal{F}_t^N$ for any $t\in\R$, where
$\mathcal{F}^N:=\{\mathcal{F}_t^N\}_{t\in\R}$ is the natural
filtration of the point process $N$,
that is,
\[
\mathcal{F}_t^{N}:=\sigma\bigl\{N(A): A\in\mathcal{B}(
\R), A\subseteq(-\infty,t]\bigr\}.
\]
Let $\{\lambda(t)\}_{t\in\R}$ be a nonnegative stochastic process
defined on $(\Omega,\mathcal{A},P)$ which is $\mathcal{F}$-adapted,
that is,
$\lambda(t)$ is $\mathcal{F}_t$-measurable for any $t\in\R$, and
such that
\[
\int_a^b\lambda(t) \,\mrr t<\infty,\qquad
\mbox{ a.s., for all }a,b\in\R.
\]
We call $\{\lambda(t)\}_{t\in\R}$ $\mathcal{F}$-stochastic
intensity of $N$ if, for any $a,b\in\R$,
\[
\mathrm{E}\bigl[N\bigl((a,b]\bigr) \mid\mathcal{F}_a]=
\mathrm{E} \biggl[\int_a^b\lambda(t)
\,\mrr t \Big| \mathcal{F}_a \biggr], \qquad\mbox{a.s.}
\]
Since one usually considers predictable stochastic intensities, we
define the predictable $\sigma$-field. Given a filtration $\mathcal
{G}:=\{\mathcal{G}_t\}_{t\in\R}\subset\mathcal{A}$,
we define the $\sigma$-field $\mathcal{P}(\mathcal{G})$ on $\R
\times\Omega$ by
\[
\mathcal{P}(\mathcal{G}):=\sigma\bigl\{(a,b]\times A: a,b\in\R, A\in
\mathcal{G}_a\bigr\}.
\]
We call $\mathcal{P}(\mathcal{G})$ predictable $\sigma$-field and
say that a real-valued stochastic process $\{X(t)\}_{t\in\R}$
is $\mathcal{G}$-predictable if the mapping $X:\R\times\Omega\to\R
$ is $\mathcal{P}(\mathcal{G})$-measurable.
A~typical $\mathcal{G}$-predictable process is a $\mathcal
{G}$-adapted process with left-continuous trajectories.

\subsection{Point processes constructed by embedding in a bivariate Poisson process}\label{subsec:emb}

Hereafter, $\overline N$ denotes a Poisson process on $\R\times\R
_+$, defined on a probability
space $(\Omega,\mathcal{A},P)$, with mean measure $\mathrm
{d}t\,\mrr z$. Let $\mathcal{F}^{\overline{N}}:=\{\mathcal
{F}_t^{\overline{N}}\}_{t\in\R}$ be the natural
filtration of $\overline{N}$, that is,
\[
\mathcal{F}_t^{\overline{N}}:=\sigma\bigl\{\overline{N}(A\times B):
A\in\mathcal{B}(\R), B\in\mathcal{B}(\R_+), A\subseteq(-\infty,t]\bigr\}.
\]
Point processes with stochastic intensity may be constructed by
embedding in a bivariate Poisson process as follows.
\begin{Lemma}\label{le:PoissEmb}
Let $f,g:\R\times\Omega\to\R_+$ be two nonnegative, $\mathcal
{P}(\mathcal{F}^{\overline{N}})$-measurable mappings
such that
\[
\int_a^b \bigl\llvert f(t)-g(t)\bigr\rrvert
\,\mrr t<\infty,\qquad\mbox{a.s., for all }a,b\in\R,
\]
set $I_t:=(\min\{f(t),g(t)\},\max\{f(t),g(t)\}]$, $t\in\R$, and
define the point process on~$\R$
\[
N(\mrr t):=\overline{N}(\mrr t\times I_t),\qquad t\in\R.
\]
Then $N$ has $\mathcal{F}^{\overline N}$-stochastic intensity $\{
\llvert  f(t)-g(t)\rrvert  \}_{t\in\R}$.
\end{Lemma}
This result is an extension of the method proposed in \cite{lewis} for
the simulation
of nonhomogeneous Poisson processes and was used, for example, in \cite
{bremaudmassoulie1}
and \cite{massoulie} to study the stability of various classes of
point processes,
including Hawkes processes.

Throughout this paper, we consider point processes $N$ on $\R$ defined by
\begin{equation}
\label{eq:stoch1} N(\mrr t):=\overline{N}\bigl(\mrr t\times\bigl(0,\lambda(t)\bigr]
\bigr),
\end{equation}
where
$\{\lambda(t)\}_{t\in\R}$ is a nonnegative process of the form
\begin{equation}
\label{eq:stoch2} \lambda(t):=\varphi(t,\overline{N}\mid_{(-\infty,t)})
\end{equation}
such that
\begin{equation}
\label{eq:stoch3} \int_a^b\lambda(s) \,\mrr s<
\infty,\qquad\mbox{a.s., for all }a,b\in\R.
\end{equation}
Here, $\varphi:\R\times\mathcal N\to\R_+$ is a measurable
functional, $\mathcal N$ denotes the space of simple and locally finite
counting measures on $\R\times\R_+$ endowed with the vague topology
(see, e.g., \cite{daley2}) and, for simplicity, with a little abuse of
notation, we denote by $\overline{N}\mid_{(-\infty,t)}$ the restriction
of $\overline{N}$ to $(-\infty,t)\times\R_+$,
that is,
\[
\overline{N}\mid_{(-\infty,t)}(A):=\overline{N}\bigl(A\cap\bigl((-\infty,t)
\times\R_+\bigr)\bigr),\qquad A\in\mathcal{B}(\R)\otimes\mathcal{B}(\R_+).
\]
Since the process $\{\overline{N}\mid_{(-\infty,t)}(A)\}_{t\in\R}$
is $\mathcal{F}^{\overline N}$-adapted and left-continuous the mapping
\[
(t,\omega)\to\overline{N}(\omega)\mid_{(-\infty,t)}(A)
\]
is $\mathcal{P}(\mathcal{F}^{\overline{N}})$-measurable for any
fixed $A\subseteq\mathcal{B}(\R)\otimes\mathcal{B}(\R_+)$.
Therefore, $\{\lambda(t)\}_{t\in\R}$ is $\mathcal{F}^{\overline
{N}}$-predictable
(see, e.g., Remark 1 in \cite{massoulie}). Consequently, by Lemma \ref
{le:PoissEmb} we deduce that
$N$ defined by \eqref{eq:stoch1}, \eqref{eq:stoch2} and \eqref
{eq:stoch3} has
$\mathcal{F}^{\overline N}$-stochastic intensity $\{\lambda(t)\}
_{t\in\R}$.

As we shall see more in detail later on, Hawkes processes may be
constructed by embedding in a bivariate Poisson process;
see \cite{bremaudmassoulie1}.

\subsection{The finite difference operator on the Poisson space and a Clark--Ocone type representation formula}

Given a measurable functional $\psi:\mathcal N\to\R$,
we define the finite difference operator $D$ by
\[
D_{(t,z)}\psi(\overline{N}):=\psi(\overline{N}+\varepsilon
_{(t,z)})-\psi(\overline{N}),
\]
where $\varepsilon_{(t,z)}$ denotes the Dirac measure at $(t,z)\in\R
\times\R_+$. We also define the $\sigma$-field
\[
\mathcal{F}_{t^-}^{\overline{N}}:=\sigma\bigl\{\overline{N}(A\times B):
A\in\mathcal{B}(\R), B\in\mathcal{B}(\R_+), A\subseteq (-\infty,t)\bigr\},\qquad t\in
\R.
\]
The following Clark--Ocone type representation formula holds; see
Theorem 1.1 in \cite{last} (see also Lemma 1.3 in \cite{wu}).
\begin{Lemma}\label{le:wu}
For any measurable functional $\psi:\mathcal N\to\R$ such that $\psi
(\overline{N})\in L^2(\Omega,\mrr P)$, we have
%
\[
\psi(\overline{N})-\mathrm{E}\bigl[\psi(\overline{N})\bigr]=\int
_{\R\times
\R_+} \mathrm{E}\bigl[D_{(t,z)}\psi(\overline{N}) \mid
\mathcal {F}_{t^-}^{\overline{N}}\bigr]\bigl(\overline{N}(\mrr t
\times\mathrm {d}z)-\mrr t\,\mrr z\bigr).
\]
\end{Lemma}
As pointed out in \cite{last} and \cite{wu}, we can (and we will)
work with a $\mathcal{P}(\mathcal{F}^{\overline{N}})\otimes\mathcal
{B}(\R_+)$-measurable version of the conditional expectation
$\mathrm{E}[D_{(t,z)}\psi(\overline{N}) \mid \mathcal
{F}_{t^-}^{\overline{N}}]$.

\section{Gaussian approximation of the first chaos of point processes with stochastic intensity}\label{sec:3}

In this section, we provide a bound for the Wasserstein distance
between a standard normal random variable $Z$
and the first chaos
\begin{equation}
\label{eq:1stchaos} \delta(u):=\int_{\R}u(t) \bigl(N(\mrr t)-
\lambda(t) \,\mrr t\bigr),
\end{equation}
being $u:\R\to\R$ a measurable function and $N$ defined by \eqref
{eq:stoch1},
\eqref{eq:stoch2} and \eqref{eq:stoch3}. We recall that, given two
random variables $X$, $Y$ defined on the same probability space, the
Wasserstein distance between $X$ and $Y$ is
\[
d_W(X,Y):=\sup_{h\in\operatorname{Lip}(1)}\bigl\llvert \mathrm{E}
\bigl[h(X)\bigr]-\mathrm{E}\bigl[h(Y)\bigr]\bigr\rrvert,
\]
where $\operatorname{Lip}(1)$ denotes the class of real-valued Lipschitz
functions with Lipschitz constant less than or equal to $1$.
We also recall that the topology induced by $d_W$ on the class of
probability measures over $\R$ is finer than the topology of weak
convergence (see, e.g., \cite{dudley}).

Following \cite{peccati}, we give a general bound for $d_W(X,Z)$,
where $X$
is an integrable random variable. Given $h\in\operatorname{Lip}(1)$, it
turns out that there exists a
twice differentiable function $f_h:\R\to\R$ so that
\begin{equation}
\label{eq:stein} h(x)-\mathrm{E}\bigl[h(Z)\bigr]=f_h'(x)-x
f_h(x),\qquad x\in\R.
\end{equation}
For a function $g:\R\to\R$, we define
$\llVert  g\rrVert  _\infty:=\sup_{x\in\R}\llvert  g(x)\rrvert  $. Equation \eqref{eq:stein} is
called Stein's equation \cite{stein} and the function $f_h$ has the following
properties:
\[
\llVert f_h\rrVert _\infty\leq2\bigl\llVert h'
\bigr\rrVert _\infty, \qquad\bigl\llVert f_h'
\bigr\rrVert _\infty\leq\sqrt{2/\pi}\bigl\llVert h'\bigr
\rrVert _\infty,\qquad\llVert f_h\rrVert _\infty\leq
2\bigl\llVert h'\bigr\rrVert _\infty;
\]
see \cite{chen}, Lemma 2.4. Since $\llVert  h'\rrVert  _\infty\leq1$ (indeed $h$ has
Lipschitz constant less than or equal to $1$), letting
$\mathcal{F}_W$ denote the class of twice differentiable functions
$f$ so that $\llVert  f\rrVert  _\infty\leq2$, $\llVert  f'\rrVert  _\infty\leq\sqrt{2/\pi}$
and $\llVert  f''\rrVert  _\infty\leq2$, we have
\begin{equation}
\label{eq:fundineq} d_W(X,Z)\leq\sup_{f\in\mathcal{F}_W}\bigl\llvert
\mathrm{E}\bigl[f'(X)-Xf(X)\bigr]\bigr\rrvert.
\end{equation}
Note that the right-hand side of \eqref{eq:fundineq} is finite since the
functions $f,f'$ are bounded and $X$ is integrable.

The following upper bound extends Corollary 3.4 in \cite{peccati} to a
class of not necessarily Poisson processes.
\begin{Theorem}\label{thm:normalskor}
Let $u:\R\to\R$ be a measurable function such that
\begin{eqnarray}
\label{eq:L1} \mathrm{E} \biggl[\int_{\R}\bigl\llvert u(t)
\bigr\rrvert \lambda(t) \,\mrr t \biggr]&<&\infty,
\\
\label{eq:L2} \mathrm{E} \biggl[\int_{\R}\bigl\llvert u(t)
\bigr\rrvert ^2\lambda(t) \,\mrr t \biggr]&<&\infty,
\\
\label{eq:L3} \int_{\R\times\R_+} \biggl(\int_t^\infty
\bigl\llvert u(s)\bigr\rrvert \mathrm {E}\bigl[\bigl\llvert D_{(t,z)}
\lambda(s)\bigr\rrvert \bigr] \,\mrr s \biggr) \,\mrr t\,\mathrm {d}z&<&\infty,
\\
\label{eq:L4} \int_{\R\times\R_+} \biggl(\int_t^\infty
\bigl\llvert u(s)\bigr\rrvert ^2\mathrm {E}\bigl[\bigl\llvert
D_{(t,z)}\lambda(s)\bigr\rrvert \bigr] \,\mrr s \biggr) \,\mrr t\,
\mathrm {d}z&<&\infty
\end{eqnarray}
and
\begin{equation}
\label{eq:L5} \int_{\R\times\R_+}\bigl\llvert u(t)\bigr\rrvert
^2 \biggl(\int_t^\infty\bigl\llvert
u(s)\bigr\rrvert \mathrm {E}\bigl[\ind_{(0,\lambda(t)]}(z)\bigl\llvert
D_{(t,z)}\lambda(s)\bigr\rrvert \bigr] \,\mathrm {d}s \biggr) \,\mrr t
\,\mrr z<\infty.
\end{equation}
In addition, assume that, for $\mrr x\,\mrr y$-almost all
$(t,z)\in\R\times\R_+$, the random function $\llvert  D_{(t,z)}\lambda
(\cdot)\rrvert  $ is a.s. locally integrable
on $(t,\infty)$ with respect to the Lebesgue measure.
Then
\begin{eqnarray}\label{eq:boundgauss}
&& d_W\bigl(\delta(u),Z\bigr)\nonumber
\\
&&\qquad \leq \sqrt{2/\pi}\mathrm{E} \biggl[\biggl
\llvert 1-\int_{\R}\bigl\llvert u(t)\bigr\rrvert
^2\lambda(t) \,\mrr t\biggr\rrvert \biggr] +\mathrm{E} \biggl[\int
_{\R}\bigl\llvert u(t)\bigr\rrvert ^3\lambda(t)
\,\mrr t \biggr]
\nonumber\\[-8pt]\\[-8pt]\nonumber
&&\quad\qquad {}+2\sqrt{2/\pi}\int_{\R\times\R_+}\bigl\llvert u(t)\bigr\rrvert
\biggl(\int_{t}^{+\infty} \bigl\llvert u(s)\bigr\rrvert
\mathrm{E}\bigl[\ind_{(0,\lambda(t)]}(z)\bigl\llvert D_{(t,z)}\lambda(s)\bigr
\rrvert \bigr] \,\mrr s \biggr) \,\mrr t\,\mrr z
\\
&&\qquad\quad{}+\int_{\R\times\R_+}\bigl\llvert u(t)\bigr
\rrvert \biggl(\int_{t}^{+\infty} \bigl\llvert u(s)\bigr
\rrvert ^2\mathrm{E}\bigl[\ind_{(0,\lambda(t)]}(z)\bigl\llvert
D_{(t,z)}\lambda (s)\bigr\rrvert \bigr]\,\mrr s \biggr) \,\mrr t
\,\mrr z,\nonumber
\end{eqnarray}
where $\delta(u)$ is defined by \eqref{eq:1stchaos}.
\end{Theorem}

\begin{Remark}
Note that if the function $u$ is bounded, then conditions \eqref
{eq:L1} and \eqref{eq:L3} imply
\eqref{eq:L2}, \eqref{eq:L4} and \eqref{eq:L5}.
\end{Remark}

\begin{pf*}{Proof of Theorem  \ref{thm:normalskor}}
We may assume
\begin{eqnarray}
\label{eq:momterzo} \mathrm{E} \biggl[\int_{\R}\bigl\llvert u(t)
\bigr\rrvert ^3\lambda(t) \,\mrr t \biggr]&<&\infty,
\\
\label{eq:infty1} \qquad\int_{\R\times\R_+}\bigl\llvert u(t)\bigr\rrvert \biggl(
\int_{t}^{+\infty} \bigl\llvert u(s)\bigr\rrvert
\mathrm{E}\bigl[\ind_{(0,\lambda(t)]}(z)\bigl\llvert D_{(t,z)}\lambda(s)\bigr
\rrvert \bigr] \,\mrr s \biggr) \,\mrr t\,\mrr z&<&\infty
\end{eqnarray}
and
\begin{equation}
\label{eq:infty2} \qquad\int_{\R\times\R_+}\bigl\llvert u(t)\bigr\rrvert \biggl(
\int_{t}^{+\infty} \bigl\llvert u(s)\bigr\rrvert
^2\mathrm{E}\bigl[\ind_{(0,\lambda(t)]}(z)\bigl\llvert D_{(t,z)}
\lambda(s)\bigr\rrvert \bigr] \,\mrr s \biggr) \,\mrr t\,\mrr z<\infty.
\end{equation}
Indeed, if one of the above terms is equal to infinity, then the claim
is trivially true.
We have
\begin{eqnarray}
\delta(u)=\int_{\R}u(t) \bigl(N(\mrr t)-\lambda(t)
\,\mrr t\bigr) &=&\int_{\R}u(t) \bigl(\overline{N}\bigl(\mrr t\times\bigl(0,\lambda (t)\bigr]\bigr)-\lambda(t) \,\mrr t\bigr)
\nonumber
\\
&=&\int_{\R\times\R_+}u(t)\ind_{(0,\lambda(t)]}(z) \bigl(\overline {N}(\mrr t\times \mrr z)-\mrr t\,\mathrm {d}z\bigr).
\nonumber
\end{eqnarray}
For any $f\in\mathcal{F}_W$, we have $f(\delta(u))\in L^2(\Omega,\mrr P)$ since $f$ is bounded.
So by Lemma~\ref{le:wu} we deduce
\[
f\bigl(\delta(u)\bigr)-\mathrm{E}\bigl[f\bigl(\delta(u)\bigr)\bigr]= \int
_{\R\times\R_+}\mathrm{E}\bigl[D_{(t,z)}f\bigl(\delta(u)\bigr)
\mid \mathcal {F}_{t^-}^{\overline{N}}\bigr]\bigl(\overline{N}(\mrr t\times\mathrm {d}z)-\mrr t\,\mrr z\bigr).
\]
For ease of notation, we set
\[
g_1(t,\omega,z):=u(t)\ind_{(0,\lambda(t,\omega)]}(z)\quad\mbox {and}\quad
g_2(t,\omega,z):=\mathrm{E}\bigl[D_{(t,z)}f\bigl(\delta(u)
\bigr) \mid \mathcal {F}_{t^-}^{\overline{N}}\bigr](\omega).
\]
By the arguments at the end of Section~\ref{subsec:emb} and the
comment after the statement of Lemma~\ref{le:wu},
we have that $g_1$ and $g_2$ are $\mathcal{P}(\mathcal{F}^{\overline
{N}})\otimes\mathcal{B}(\R_+)$-measurable.
Note that, due to assumptions \eqref{eq:L1} and \eqref{eq:L2}, $g_1$
is integrable and square integrable with respect to $\mathrm
{d}t\,\mrr z\,\mrr P(\omega)$.
We shall check later on that
\begin{equation}
\label{eq:albarb2} \qquad g_2\mbox{ is integrable and square integrable with
respect to }\mathrm {d}t\,\mrr z\,\mrr P(\omega).
\end{equation}
So by Theorem 3 in \cite{bremaudmassoulie} [formulas (19)~and~(20)], we have
\begin{eqnarray}
\mathrm{E}\bigl[\delta(u)f\bigl(\delta(u)\bigr)\bigr]&=&\mathrm{E}\bigl[\delta(u)
\bigl(f\bigl(\delta (u)\bigr)-\mathrm{E}\bigl[f\bigl(\delta(u)\bigr)\bigr]\bigr)
\bigr]
\nonumber
\\
&=&\mathrm{E} \biggl[ \biggl(\int_{\R\times\R_+}g_1(t,z)
\bigl(\overline {N}(\mrr t\times \mrr z)-\mrr t\,\mrr z\bigr)
\biggr)
\nonumber
\\
&&{} \times \biggl(\int_{\R\times\R_+}g_2(t,z) \bigl(
\overline{N}(\mathrm {d}t\times \mrr z)-\mrr t\,\mrr z\bigr) \biggr)
\biggr]
\nonumber
\\
&=&\mathrm{E} \biggl[\int_{\R\times\R_+}g_1(t,z)g_2(t,z)\,
\mathrm {d}t\,\mrr z \biggr].
\nonumber
\end{eqnarray}
By the Taylor formula, we deduce
\begin{eqnarray}
D_{(t,z)}f\bigl(\delta(u)\bigr)&=&f\bigl(\delta(u)+D_{(t,z)}
\delta(u)\bigr)-f\bigl(\delta (u)\bigr)\label{eq:findiff}
\\
%
&=&f'
\bigl(\delta(u)\bigr)D_{(t,z)}\delta(u)+R\bigl(D_{(t,z)}\delta(u)
\bigr),\label{eq:taylorsecond}
\end{eqnarray}
where the rest $R$ satisfies $\llvert  R(y)\rrvert  \leq y^2$ since $\llVert  f''\rrVert  _\infty
\leq2$. Since $f$ is bounded, by~\eqref{eq:findiff} we have that
$g_2$ is a.s. bounded,
and so by the standard properties of the conditional expectation [note
that $\lambda(t)$ is $\mathcal{F}_{t^-}^{\overline N}$-measurable]
and Fubini's theorem we deduce
\begin{eqnarray}
\mathrm{E}\bigl[\delta(u)f\bigl(\delta(u)\bigr)\bigr] =\mathrm{E} \biggl[\int
_{\R\times\R_+}g_1(t,z)D_{(t,z)}f\bigl(\delta (u)
\bigr) \,\mrr t\,\mrr z \biggr].
\nonumber
\end{eqnarray}
Consequently, by \eqref{eq:taylorsecond},
%
\begin{eqnarray*}
&&\bigl\llvert \mathrm{E}\bigl[f'\bigl(\delta(u)\bigr)-\delta(u)f
\bigl(\delta(u)\bigr)\bigr]\bigr\rrvert
\\
&&\qquad =\biggl\llvert \mathrm {E}
\biggl[f'\bigl(\delta(u)\bigr)-\int_{\R\times\R
_+}g_1(t,z)D_{(t,z)}f
\bigl(\delta(u)\bigr) \,\mrr t\,\mrr z \biggr]\biggr\rrvert
\nonumber
\\
&&\qquad \leq\biggl\llvert \mathrm{E}
\biggl[f'\bigl(\delta(u)\bigr) \biggl(1-\int_{\R\times
\R_+}g_1(t,z)D_{(t,z)}
\delta(u) \,\mrr t\,\mrr z \biggr) \biggr]\biggr\rrvert 
\\
&&\quad\qquad{} +
\biggl\llvert \mathrm{E} \biggl[\int_{\R\times\R
_+}g_1(t,z)R
\bigl(D_{(t,z)}\delta(u)\bigr) \,\mrr t\,\mrr z \biggr]\biggr\rrvert
\nonumber
\\
&&\quad\qquad {} \leq\sqrt{2/\pi}\mathrm{E} \biggl[\biggl\llvert 1-\int_{\R\times\R
_+}g_1(t,z)D_{(t,z)}
\delta(u) \,\mrr t\,\mrr z\biggr\rrvert \biggr]
\\
&&\quad\qquad{}+ \mathrm{E} \biggl[\int
_{\R\times\R_+}\bigl\llvert g_1(t,z)\bigr\rrvert \bigl
\llvert D_{(t,z)}\delta (u)\bigr\rrvert ^2 \,\mrr t
\,\mrr z \biggr].
\nonumber
\end{eqnarray*}
Therefore, using the basic inequality \eqref{eq:fundineq}, we have
\begin{eqnarray}
d_W\bigl(\delta(u),Z\bigr) &\leq&\sqrt{2/\pi}\mathrm{E} \biggl[
\biggl\llvert 1-\int_{\R\times\R
_+}g_1(t,z)D_{(t,z)}
\delta(u) \,\mrr t\,\mrr z\biggr\rrvert \biggr]
\nonumber
\\
&&{} +\mathrm{E} \biggl[\int_{\R\times\R_+}\bigl\llvert
g_1(t,z)\bigr\rrvert \bigl\llvert D_{(t,z)}\delta (u)\bigr
\rrvert ^2 \,\mrr t\,\mrr z \biggr].
\nonumber
\end{eqnarray}
We shall check later on that
\begin{equation}
\label{eq:Fub1} \mathrm{E} \biggl[\int_{\R\times\R_+}\bigl\llvert
g_1(t,z)\bigr\rrvert \bigl\llvert D_{(t,z)}\delta (u)\bigr
\rrvert \,\mrr t\,\mrr z \biggr]<\infty
\end{equation}
and
\begin{equation}
\label{eq:Fub2} \mathrm{E} \biggl[\int_{\R\times\R_+}\bigl\llvert
g_1(t,z)\bigr\rrvert \bigl\llvert D_{(t,z)}\delta (u)\bigr
\rrvert ^2 \,\mrr t\,\mrr z \biggr]<\infty.
\end{equation}
So the above upper bound on $d_W(\delta(u),Z)$ is nontrivial. For
$\mrr x\,\mrr y$-almost all $(t,z)\in\R\times\R_+$, we have
\begin{eqnarray}
D_{(t,z)}\delta(u) 
&=&D_{(t,z)}
\biggl(\int_{\R\times\R_+}g_1(s,v)\overline{N}(\mathrm {d}s
\times \mrr v) \biggr) -D_{(t,z)} \biggl(\int_{\R\times\R_+}g_1(s,v)
\,\mrr s\,\mathrm {d}v \biggr).
\nonumber
\end{eqnarray}
Computing separately these two finite differences and writing $\varphi
_t(\overline{N}\mid_{(-\infty,t)})$ in place of $\varphi(t,\overline
{N}\mid_{(-\infty,t)})$
for ease of notation, we have
\begin{eqnarray*}
&& D_{(t,z)} \biggl(\int_{\R\times\R_+}g_1(s,v)
\overline {N}(\mrr s\times \mrr v) \biggr)
\\
&&\qquad =  D_{(t,z)} \biggl(\int_{\R\times\R_+}\ind_{s\leq t}
u(s)\ind_{(0,\varphi_s(\overline{N}\mid_{(-\infty,s)})]}(v)\overline {N}(\mrr s\times \mrr v)
\\
&&\quad\qquad{} + \int_{\R\times\R_+}\ind_{s>t} u(s)\ind_{(0,\varphi_s(\overline{N}\mid_{(-\infty,s)})]}(v)
\overline {N}(\mrr s\times \mrr v) \biggr)
\\
&&\qquad  =  \int_{\R\times\R_+}\ind_{s\leq t} u(s)\ind_{(0,\varphi_s((\overline{N}+\varepsilon
_{(t,z)})\mid_{(-\infty,s)})]}(v)
(\overline{N}+\varepsilon _{(t,z)}) (\mrr s\times \mrr v)
\\
&&\quad\qquad{} + \int_{\R\times\R_+}\ind_{s>t} u(s)\ind_{(0,\varphi_s((\overline{N}+\varepsilon
_{(t,z)})\mid_{(-\infty,s)})]}(v)
(\overline{N}+\varepsilon _{(t,z)}) (\mrr s\times \mrr v)
\\
&&\quad\qquad{} -\int_{\R\times\R_+}\ind_{s\leq t} u(s)\ind_{(0,\varphi_s(\overline{N}\mid_{(-\infty,s)})]}(v)
\overline {N}(\mrr s\times \mrr v)
\\
&&\quad\qquad{} - \int_{\R\times\R_+}\ind_{s>t} u(s)\ind_{(0,\varphi_s(\overline{N}\mid_{(-\infty,s)})]}(v)
\overline {N}(\mrr s\times \mrr v)
\nonumber
\\
&&\qquad =  \int_{\R\times\R_+}\ind_{s\leq t} u(s)\ind_{(0,\varphi_s(\overline{N}\mid_{(-\infty,s)})]}(v)
(\overline {N}+\varepsilon_{(t,z)}) (\mrr s\times \mrr v)
\\
&&\quad\qquad{} + \int_{\R\times\R_+}\ind_{s>t} u(s)\ind_{(0,\varphi_s(\overline{N}\mid_{(-\infty,s)}+\varepsilon
_{(t,z)})]}(v)
\overline{N}(\mrr s\times \mrr v)
\\
&&\quad\qquad{} -\int_{\R\times\R_+}\ind_{s\leq t} u(s)\ind_{(0,\varphi_s(\overline{N}\mid_{(-\infty,s)})]}(v)
\overline {N}(\mrr s\times \mrr v)
\\
&&\quad\qquad{} -\int_{\R\times\R_+}\ind_{s>t} u(s)\ind_{(0,\varphi_s(\overline{N}\mid_{(-\infty,s)})]}(v)
\overline {N}(\mrr s\times \mrr v)
\\
&&\qquad =  g_1(t,z)+\int_{(t,\infty)}u_{(t,z)}(s)N_{(t,z)}(\mrr s),
\end{eqnarray*}
where for $s>t$
\begin{eqnarray*}
u_{(t,z)}(s) &:=&\operatorname{sign}\bigl(\varphi_s(\overline{N}
\mid_{(-\infty ,s)}+\varepsilon_{(t,z)})-\varphi_s(
\overline{N}\mid_{(-\infty,s)})\bigr)u(s)
\\
&=&\operatorname{sign}\bigl(D_{(t,z)}
\lambda(s)\bigr)u(s),
\\
N_{(t,z)}(\mrr s)&:=&\overline{N}\bigl(\mrr s\times \bigl(
\varphi_s(\overline{N}\mid_{(-\infty,s)}+\varepsilon_{(t,z)})
\wedge \varphi_s(\overline{N}\mid_{(-\infty,s)}),
\\
&&{} \varphi_s(\overline{N}\mid_{(-\infty,s)}+\varepsilon_{(t,z)})
\vee \varphi_s(\overline{N}\mid_{(-\infty,s)})\bigr]\bigr).
\nonumber
\end{eqnarray*}
Here, for ease of notation, we denoted by $a\wedge b$ and $a\vee b$ the
minimum and the maximum between $a,b\in\R$, respectively.
Moreover,
\begin{eqnarray*}
D_{(t,z)} \biggl(\int_{\R\times\R_+} g_1(s,v)
\,\mrr s\,\mrr v \biggr) &=&\int
_{t}^{+\infty} u(s)D_{(t,z)}\lambda(s)
\,\mrr s 
\\
&=&\int_{t}^{+\infty}
u_{(t,z)}(s)\bigl\llvert D_{(t,z)}\lambda(s)\bigr\rrvert
\,\mrr s.
\nonumber
\end{eqnarray*}
Therefore,
\begin{equation}
\label{eq:Dtz} D_{(t,z)}\delta(u)=g_1(t,z)+
\delta_{(t,z)}(u),
\end{equation}
where
\begin{eqnarray}
\delta_{(t,z)}(u)&:=&\int_{(t,\infty)}u_{(t,z)}(s)
\bigl(N_{(t,z)}(\mrr s) 
%
-
\bigl\llvert D_{(t,z)}\lambda(s)\bigr\rrvert \,\mrr s\bigr).
\nonumber
\end{eqnarray}
Combining \eqref{eq:Dtz} with the previous bound on $d_W(\delta
(u),Z)$, we deduce
\begin{eqnarray}\label{eq:diseq1}
d_W\bigl(\delta(u),Z\bigr)&\leq&\sqrt{2/\pi}\mathrm{E} \biggl[\biggl
\llvert 1-\int_{\R}\bigl\llvert u(t)\bigr\rrvert
^2\lambda(t) \,\mrr t\biggr\rrvert \biggr]\nonumber
\\
&&{} + \sqrt{2/\pi}\mathrm{E} \biggl[\int_{\R\times\R_+}\bigl\llvert u(t)\bigr
\rrvert \ind _{(0,\lambda(t)]}(z)\bigl\llvert \delta_{(t,z)}(u)\bigr\rrvert
\,\mrr t\,\mathrm {d}z \biggr]
\nonumber
\\
&&{}+\mathrm{E} \biggl[\int_{\R}\bigl\llvert u(t)\bigr\rrvert
^3\lambda(t) \,\mrr t \biggr] 
\\
&&{}
+2\mathrm{E}
\biggl[\int_{\R\times\R_+}\bigl\llvert u(t)\bigr\rrvert u(t)
\ind_{(0,\lambda(t)]}(z) \delta_{(t,z)}(u) \,\mrr t\,\mrr z \biggr]
\nonumber
\\
&&{}+\mathrm{E} \biggl[\int_{\R\times\R_+}\bigl\llvert u(t)\bigr\rrvert
\ind_{(0,\lambda(t)]}(z) \bigl\llvert \delta_{(t,z)}(u)\bigr\rrvert
^2 \,\mrr t\,\mrr z \biggr].\nonumber
\end{eqnarray}
We shall check later on that
\begin{equation}
\label{eq:torsa1} \mathrm{E} \biggl[\int_{\R\times\R_+}\bigl\llvert
g_1(t,z)\bigr\rrvert ^2 \bigl\llvert D_{(t,z)}
\delta(u)\bigr\rrvert \,\mrr t\,\mrr z \biggr]<\infty,
\end{equation}
and so by \eqref{eq:L2}, \eqref{eq:momterzo}, \eqref{eq:Fub1},
\eqref{eq:Fub2},
\eqref{eq:Dtz} and \eqref{eq:torsa1} the bound \eqref{eq:diseq1} is
nontrivial. By Lemma~\ref{le:PoissEmb},
for $\mrr x\,\mrr y$-almost all $(t,z)\in\R\times\R_{+}$,
the point process $N_{(t,z)}$ on $(t,\infty)$
has $\{\mathcal{F}_s^{\overline{N}}\}_{s>t}$-stochastic intensity $\{
\llvert  D_{(t,z)}\lambda(s)\rrvert  \}_{s>t}$. Indeed, the mapping
\[
(t,\infty)\times\Omega\ni(s,\omega)\mapsto\bigl\llvert D_{(t,z)}\lambda
(s,\omega)\bigr\rrvert \in\R
\]
is $\mathcal{P}(\{\mathcal{F}_s^{\overline{N}}\}_{s>t})$-measurable
and (by assumption), for $P$-almost all
$\omega$, it is locally integrable in $s$ with respect to the Lebesgue
measure. We note that
\begin{eqnarray}
\ind_{(0,\lambda(t)]}(z) \delta_{(t,z)}(u) 
=
\int_{(t,\infty)}\ind_{(0,\lambda
(t)]}(z)u_{(t,z)}(s)
\bigl(N_{(t,z)}(\mrr s) 
%
%
-
\bigl\llvert D_{(t,z)}\lambda(s)\bigr\rrvert \,\mrr s\bigr)
\nonumber
\end{eqnarray}
and the mapping
\[
(t,\infty)\times\Omega\ni(s,\omega)\to\ind_{(0,\lambda(t,\omega
)]}(z)u_{(t,z)}(s,
\omega)
\]
is $\mathcal{P}(\{\mathcal{F}_{s}^{\overline{N}}\}
_{s>t})$-measurable. By \eqref{eq:L3} and \eqref{eq:L4}, we have
\begin{eqnarray}
\int_{t}^{+\infty} \bigl\llvert u(s)\bigr\rrvert
\mathrm{E}\bigl[\ind_{(0,\lambda(t)]}(z)\bigl\llvert D_{(t,z)}\lambda(s)\bigr
\rrvert \bigr] \,\mrr s<\infty,
\nonumber
\\
\eqntext{\mbox{for }\mrr x\,\mrr y\mbox
{-almost all }(t,z)\in\R\times\R_+}
\end{eqnarray}
and
\begin{eqnarray}
\int_{t}^{+\infty} \bigl\llvert u(s)\bigr\rrvert
^2\mathrm{E}\bigl[\ind_{(0,\lambda(t)]}(z)\bigl\llvert D_{(t,z)}
\lambda(s)\bigr\rrvert \bigr] \,\mrr s<\infty,\nonumber
\\
\eqntext{\mbox{for }\mrr x
\,\mrr y\mbox {-almost all }(t,z)\in\R\times\R_+.}
\end{eqnarray}
Therefore, by Theorem 3 in \cite{bremaudmassoulie} [formulas (19)~and~(20)], for $\mrr x\,\mrr y$-almost all $(t,z)\in\R
\times\R_+$, we have
\begin{eqnarray}
\mathrm{E}\bigl[\ind_{(0,\lambda(t)]}(z) \delta_{(t,z)}(u)
\bigr]=0\label{eq:dirfc1}
\end{eqnarray}
and
\begin{eqnarray}
\qquad \mathrm{E}\bigl[\ind_{(0,\lambda(t)]}(z) \bigl\llvert \delta_{(t,z)}(u)
\bigr\rrvert ^2\bigr] 
&=&\int
_{t}^{\infty}\bigl\llvert u(s)\bigr\rrvert
^2\mathrm{E}\bigl[\ind_{(0,\lambda
(t)]}(z)\bigl\llvert D_{(t,z)}
\lambda(s)\bigr\rrvert \bigr] \,\mrr s.\label{eq:dirfc2}
\end{eqnarray}
By the triangular inequality and formula (19) in \cite
{bremaudmassoulie}, we have
\begin{eqnarray}
\mathrm{E}\bigl[\ind_{(0,\lambda(t)]}(z) \bigl\llvert \delta_{(t,z)}(u)
\bigr\rrvert \bigr]&\leq& 
2\int
_{t}^{\infty}\bigl\llvert u(s)\bigr\rrvert \mathrm{E}
\bigl[\ind_{(0,\lambda
(t)]}(z)\bigl\llvert D_{(t,z)}\lambda(s)\bigr\rrvert
\bigr] \,\mrr s.\label{eq:dirfc3}
\end{eqnarray}
Inequality \eqref{eq:boundgauss} follows
combining \eqref{eq:diseq1} with \eqref{eq:dirfc1}, \eqref
{eq:dirfc2} and \eqref{eq:dirfc3}.

It remains to prove the integrability conditions \eqref{eq:albarb2},
\eqref{eq:Fub1}, \eqref{eq:Fub2} and \eqref{eq:torsa1}.
Since $f\in\mathcal{F}_W$, then it is
Lipschitz continuous with Lipschitz constant less than or equal to $1$.
Therefore,
by \eqref{eq:findiff} we have $\llvert  D_{(t,z)}f(\delta(u))\rrvert  \leq
\llvert  D_{(t,z)}\delta(u)\rrvert  $,
and so to prove \eqref{eq:albarb2} it suffices to check
\[
\int_{\R\times\R_+}\mathrm{E}\bigl[\bigl\llvert D_{(t,z)}
\delta(u)\bigr\rrvert \bigr] \,\mathrm {d}t\,\mrr z<\infty\quad\mbox{and}\quad
\int_{\R\times\R_+}\mathrm{E}\bigl[
\bigl\llvert D_{(t,z)}\delta(u)\bigr\rrvert ^2\bigr] \,\mathrm
{d}t\,\mrr z<\infty.
\]
Using relation \eqref{eq:Dtz} and formula (19) in \cite
{bremaudmassoulie}, we have
\begin{eqnarray*}
&& \int_{\R\times\R_+}\mathrm{E}\bigl[\bigl\llvert D_{(t,z)}
\delta(u)\bigr\rrvert \bigr] \,\mathrm {d}t\,\mrr z
\\
&&\qquad  \leq\mathrm{E} \biggl[\int
_{\R}\bigl\llvert u(t)\bigr\rrvert \lambda(t) \,\mrr t
\biggr]+ 2\mathrm{E} \biggl[\int_{\R\times\R_+} \biggl(\int
_t^\infty \bigl\llvert u(s)\bigr\rrvert \bigl
\llvert D_{(t,z)}\lambda(s)\bigr\rrvert \,\mrr s \biggr) \,\mrr t\,
\mathrm {d}z \biggr]
\nonumber
\end{eqnarray*}
and this latter term is finite due to assumptions \eqref{eq:L1} and
\eqref{eq:L3}. Using again relation \eqref{eq:Dtz} and formula (20)
in \cite{bremaudmassoulie}, we have
\begin{eqnarray*}
&& \int_{\R\times\R_+}\mathrm{E}\bigl[\bigl\llvert D_{(t,z)}
\delta(u)\bigr\rrvert ^2\bigr] \,\mathrm {d}t\,\mrr z
\\
&&\qquad  \leq 2
\mathrm{E} \biggl[\int_{\R}\bigl\llvert u(t)\bigr\rrvert
^2\lambda(t) \,\mathrm {d}t \biggr]
+2\mathrm{E} \biggl[\int_{\R\times\R_+} \biggl(\int
_t^\infty \bigl\llvert u(s)\bigr\rrvert
^2\bigl\llvert D_{(t,z)}\lambda(s)\bigr\rrvert \,\mrr s
\biggr) \,\mathrm {d}t\,\mrr z \biggr]
\nonumber
\end{eqnarray*}
and this latter term is finite due to assumptions \eqref{eq:L2} and
\eqref{eq:L4}. By \eqref{eq:Dtz} and \eqref{eq:dirfc3}, we have
\begin{eqnarray*}
&&\mathrm{E} \biggl[\int_{\R\times\R_+}\bigl\llvert
g_1(t,z)\bigr\rrvert \bigl\llvert D_{(t,z)}\delta (u)\bigr
\rrvert \,\mrr t\,\mrr z \biggr]
\nonumber
\\
&&\qquad \leq \mathrm{E} \biggl[\int_{\R}\bigl\llvert u(t)\bigr
\rrvert ^2\lambda(t) \,\mrr t \biggr]+\mathrm{E} \biggl[\int
_{\R\times\R_+} \bigl\llvert g_1(t,z)\bigr\rrvert \bigl
\llvert \delta_{(t,z)}(u)\bigr\rrvert \,\mrr t\,\mrr z \biggr]
\\
&&\qquad \leq\mathrm{E} \biggl[\int_{\R}\bigl\llvert u(t)\bigr
\rrvert ^2\lambda(t) \,\mathrm {d}t \biggr]\nonumber
\\
&&\quad\qquad{} +2\int_{\R\times\R_+}
\bigl\llvert u(t)\bigr\rrvert \biggl(\int_{t}^{\infty}
\bigl\llvert u(s)\bigr\rrvert \mathrm{E}\bigl[\ind_{(0,\lambda
(t)]}(z)\bigl\llvert
D_{(t,z)}\lambda(s)\bigr\rrvert \bigr] \,\mrr s \biggr) \,\mrr t
\,\mrr z,
\nonumber
\end{eqnarray*}
and \eqref{eq:Fub1} follows by \eqref{eq:L2} and \eqref{eq:infty1}.
Similarly, by \eqref{eq:Dtz} and \eqref{eq:dirfc2} we have
\begin{eqnarray*}
&&\mathrm{E} \biggl[\int_{\R\times\R_+}\bigl\llvert
g_1(t,z)\bigr\rrvert \bigl\llvert D_{(t,z)}\delta (u)\bigr
\rrvert ^2 \,\mrr t\,\mrr z \biggr]
\nonumber
\\
&&\qquad \leq 2\mathrm{E} \biggl[\int_{\R}\bigl\llvert u(t)\bigr
\rrvert ^3\lambda(t) \,\mathrm {d}t \biggr]+2\mathrm{E} \biggl[\int
_{\R\times\R_+} \bigl\llvert g_1(t,z)\bigr\rrvert \bigl
\llvert \delta_{(t,z)}(u)\bigr\rrvert ^2 \,\mrr t\,\mrr z
\biggr]
\nonumber
\\
&&\qquad =2\mathrm{E} \biggl[\int_{\R}\bigl\llvert u(t)\bigr\rrvert
^3\lambda(t) \,\mrr t \biggr]
\\
&&\quad\qquad{}+2\int_{\R\times\R_+} \bigl
\llvert u(t)\bigr\rrvert \biggl(\int_{t}^{\infty}\bigl
\llvert u(s)\bigr\rrvert ^2\mathrm{E}\bigl[\ind_{(0,\lambda
(t)]}(z)\bigl
\llvert D_{(t,z)}\lambda(s)\bigr\rrvert \bigr] \,\mrr s \biggr)
\,\mrr t\,\mrr z,
\nonumber
\end{eqnarray*}
and \eqref{eq:Fub2} follows by \eqref{eq:momterzo} and \eqref
{eq:infty2}. Finally, \eqref{eq:torsa1}
may be checked similarly to \eqref{eq:Fub1}, but using
\eqref{eq:momterzo} and \eqref{eq:L5} in place of \eqref{eq:L2} and
\eqref{eq:infty1}, respectively.
The proof is completed.
\end{pf*}

\section{Application to stationary nonlinear Hawkes processes}\label{sec:4}

A nonlinear Hawkes process with parameters $(\phi,h)$ is a point
process $N$ on $\R$ with $\mathcal{F}^N$-stochastic intensity of the form
\begin{equation}
\label{eq:hawkes} t\mapsto\phi \biggl(\int_{(-\infty,t)}h(t-s)N(\mrr s) \biggr),\qquad t\in\R,
\end{equation}
where $\phi:\R\to\R_+$ and $h:\R_+\to\R$ are measurable
functions. A particular case is the self-exciting process (or linear
Hawkes process) with parameters $(\nu,h)$,
for which $\phi(x):=\nu+x$, for some constant $\nu>0$, and $h$ is
nonnegative.

In the seminal paper \cite{bremaudmassoulie1}, the authors
proved that if $\phi$ is Lipschitz continuous with Lipschitz constant
$\alpha$ such that
$\alpha\mu<1$, where $\mu:=\llVert  h\rrVert  _{L^1(\R_+,\mrr x)}$, then
there exists a unique stationary distribution of $N$
with dynamics \eqref{eq:hawkes} and finite intensity $\lambda
:=\mathrm{E}[N((0,1])]$.
The stationary solution is constructed by embedding in a bivariate
Poisson process, as follows. Define recursively the processes
$\lambda^{(0)}\equiv0$,
\[
N^{(n)}(\mrr t):=\overline{N}\bigl(\mrr t\times\bigl(0,\lambda^{(n)}(t)\bigr]\bigr)
\]
and
\[
\lambda^{(n+1)}(t):=
\phi \biggl(\int_{(-\infty
,t)}h(t-s)N^{(n)}(\mrr s)
\biggr),
\]
$n\geq0$, $t\in\R$, where $\overline{N}$ is a Poisson process on
$\R\times\R_+$ with mean measure $\mrr t\,\mrr z$.
It turns out that, for any fixed $n\geq0$, the point process $N^{(n)}$
is stationary and
$\{\lambda^{(n)}(t)\}_{t\in\R}$ is an $\mathcal{F}^{\overline
{N}}$-stochastic intensity of $N^{(n)}$.
It is then proved that $N^{(n)}((a,b])\to N((a,b])$ and $\lambda
^{(n)}(t)\to\lambda(t)$ a.s., for any $a,b,t\in\R$, and the
limiting process is stationary and satisfies
\begin{eqnarray}
N(\mrr t)=\overline{N}\bigl(\mrr t\times\bigl(0,\lambda (t)\bigr]\bigr),\qquad
\lambda(t)=\phi \biggl(\int_{(-\infty,t)}h(t-s)N(\mrr s) \biggr),
\qquad t\in\R
\nonumber
\end{eqnarray}
and $\lambda\in(0,\infty)$. Note that $\lambda(t)=\varphi
(t,\overline{N}\mid_{(-\infty,t)})$, for some functional
$\varphi:\mathbb R\times\mathcal N\to\mathbb R_+$ satisfying
\begin{eqnarray}\label{eq:functional}
&& \varphi(t,\overline{N}\mid_{(-\infty,t)})
\nonumber\\[-8pt]\\[-8pt]\nonumber
&&\qquad :=\phi \biggl(\int
_{(-\infty
,t)\times\R_+}\ind_{(0,\varphi(s,\overline{N}\mid_{(-\infty
,s)})]}(z)h(t-s)\overline{N}(\mrr s
\times \mrr z) \biggr).
\end{eqnarray}
Then by Lemma~\ref{le:PoissEmb} it follows that $N$ is a point process
on $\R$
with $\mathcal{F}^{N}$-stochastic intensity $\{\lambda(t)\}_{t\in\R
}$. In conclusion, $N$
is a stationary nonlinear Hawkes process with parameters $(\phi,h)$
and finite intensity.

\subsection{Explicit Gaussian bound for the first chaos of nonlinear Hawkes processes}\label{subsec:4/1}

The following explicit Gaussian bound holds.

\begin{Theorem}\label{thm:nonlinGauss}
Assume that $h:\R_+\to[0,\infty)$ is locally bounded and $\phi
:[0,\infty)\to[0,\infty)$, $\phi(0)>0$, is nondecreasing and
Lipschitz continuous,
with Lipschitz constant $\alpha$ such that $\alpha\mu<1$. Let $N$ be
a stationary nonlinear Hawkes process with parameters $(\phi,h)$ and
finite intensity
$\lambda\in(0,\infty)$.
If $u\in L^1(\R,\mrr x)$, then
\begin{eqnarray}
d_W\bigl(\delta(u),Z\bigr)\leq\mathfrak{N},\label{eq:thm41explicit}
\end{eqnarray}
where
\begin{eqnarray}\label{eq:UBnonlin}
\mathfrak{N}&:=&\sqrt{2/\pi}\max \biggl\{\bigl\llvert 1-\phi(0)\llVert u\rrVert
_{L^2(\R
,\mrr x)}^2\bigr\rrvert,\biggl\llvert 1-\frac{\phi(0)}{1-\alpha\mu}
\llVert u\rrVert _{L^2(\R
,\mrr x)}^2\biggr\rrvert \biggr\}\nonumber
\\
&&{} +
\frac{\phi(0)}{1-\alpha\mu}\llVert u\rrVert _{L^3(\R,\mrr x)}^3
+\frac{2\sqrt{2/\pi}\phi(0)\alpha\mu(2-\alpha\mu)}{(1-\alpha
\mu)^2}\llVert u\rrVert _{L^2(\R,\mrr x)}^2
\\
&&{} +\frac{\phi(0)\alpha\mu}{(1-\alpha\mu)^2}\llVert u\rrVert _{L^2(\R,\mathrm
{d}x)}\bigl
\llVert u^2\bigr\rrVert _{L^2(\R,\mrr x)}.\nonumber
\end{eqnarray}
\end{Theorem}

\begin{Remark}\label{eq:remnonstat}
Suppose that $\phi$ and $h$ satisfy the assumptions of Theorem~\ref
{thm:nonlinGauss}.
One says that $N'$ is a $($nonstationary$)$ nonlinear Hawkes process on
$\R_+$ with parameters $(\phi,h)$ and initial condition $\mathcal{IC}$
$($see \cite{bremaudmassoulie1}$)$
if $N'$ has stochastic intensity
\[
\lambda'(t):=\phi \biggl(\int_{(-\infty,t)}h(t-s)N'(
\mathrm {d}s) \biggr),\qquad t>0
\]
on $\R_+$ and $N'$ satisfies the condition $\mathcal{IC}$ on $\R_-$.
If $u'\in L^1(\R_+,\mrr x)$,
following the lines of the proof of Theorem~\ref{thm:nonlinGauss}, one
can show $($without major difficulties$)$ that the
Gaussian bound \eqref{eq:thm41explicit} holds replacing $\delta(u)$ with
\[
\delta'\bigl(u'\bigr):=\int_{\R_+}u'(t)
\bigl(N'(\mrr t)-\lambda'(t) \,\mrr t\bigr)
\]
and replacing $u$ with $u'$ $($and $\R$ with $\R_+)$ in the
expression of $\mathfrak N$.
\end{Remark}

Let $N$ be a stationary nonlinear Hawkes process with parameters $(\phi
,h)$ which satisfy the assumptions of Theorem~\ref{thm:nonlinGauss}.
Since $h\geq0$ and $\phi$ is nondecreasing and Lipschitz continuous
with Lipschitz constant $\alpha$, we have
\begin{equation}
\label{eq:disstocint} \phi(0)\leq\lambda(t)\leq\phi(0)+\alpha\int_{(-\infty
,t)}h(t-s)N(\mrr s),\qquad t\in\R.
\end{equation}
Taking the mean, we deduce $\phi(0)\leq\lambda\leq\phi(0)+\lambda
\alpha\mu$, and so
\begin{equation}
\label{eq:disintens} \phi(0)\leq\lambda\leq\frac{\phi(0)}{1-\alpha\mu}.
\end{equation}
Given an integrable function $u$, one may think of approximating the quantity
\[
\int_{\R}u(t)\lambda(t) \,\mrr t
\]
with its expectation $\lambda\int_{\R}u(t) \,\mrr t$.
Unfortunately, in general the intensity $\lambda$ is not known explicitly
(unless we consider the linear case which is treated in the next
section). However, it may be estimated by, for example, Monte Carlo
simulation (using the ergodic theorem). For a fixed positive constant
$\widehat{\lambda}\in[\phi(0),\phi(0)(1-\alpha\mu)^{-1}]$, to be
interpreted as an estimate of the intensity $\lambda$, we define the
``approximated'' first chaos by
\[
\delta_a(u):=\int_{\R}u(t) \bigl(N(\mrr t)-\widehat{\lambda }\,\mrr t\bigr).
\]
The following explicit Gaussian bound holds.

\begin{Theorem}\label{thm:nonlinGaussappr}
Under assumptions and notation of Theorem~\ref{thm:nonlinGauss}, we have
\begin{eqnarray}
d_W\bigl(\delta_a(u),Z\bigr)&\leq&\mathfrak{N}
+
\frac{2\phi(0)\alpha\mu}{1-\alpha\mu}\llVert u\rrVert _{L^1(\R,\mathrm
{d}x)}.\label{eq:thm41explicitappr}
\end{eqnarray}
\end{Theorem}

\begin{Remark}\label{re:linear}
In the case of stationary linear Hawkes processes, the bounds \eqref
{eq:thm41explicit} and \eqref{eq:thm41explicitappr} may be
$($slightly$)$ improved due to the knowledge
of the intensity~$\lambda$, see Theorem~\ref{thm:linGauss1}.
Moreover, alternative bounds may be obtained by using the spectral
theory of self-exciting processes; see Theorem
\ref{thm:linGaussappvecchio}.
\end{Remark}

The proofs of Theorems~\ref{thm:nonlinGauss} and~\ref
{thm:nonlinGaussappr} are given in Section~\ref{sec:proofs1}.

\subsection{A quantitative central limit theorem for nonlinear Hawkes processes}\label{subsec:4/2}

The following quantitative central limit theorem in the Wasserstein
distance is an immediate consequence of
Theorems~\ref{thm:nonlinGauss} and~\ref{thm:nonlinGaussappr}.

\begin{Corollary}\label{cor:QCLTnonlinear}
For $\varepsilon>0$, assume that $h_\varepsilon:\R_+\to[0,\infty)$
is locally bounded and
$\phi_\varepsilon:[0,\infty)\to[0,\infty)$, $\phi_\varepsilon
(0)>0$, is nondecreasing and Lipschitz continuous,
with Lipschitz constant $\alpha_\varepsilon$ such that $\alpha
_\varepsilon\mu_\varepsilon<1$, where $\mu_\varepsilon:=\int_0^\infty h_\varepsilon(t) \,\mrr t$.
Let $N_\varepsilon$ be a stationary nonlinear Hawkes process with parameters
$(\phi_\varepsilon,h_\varepsilon)$ and finite intensity $\lambda
_\varepsilon\in\R_+$, and take $u_\varepsilon\in L^1(\R,\mathrm
{d}x)$. Then:

\begin{longlist}[(ii)]
\item[(i)]
\begin{eqnarray}
d_W\bigl(\delta^{(\varepsilon)}(u_\varepsilon),Z\bigr)\leq
\mathfrak {N}_\varepsilon,\qquad\varepsilon>0\label{eq:bound1bis}
\end{eqnarray}
and
\begin{eqnarray}
d_W\bigl(\delta_a^{(\varepsilon)}(u_\varepsilon),Z
\bigr)\leq\mathfrak {N}_\varepsilon +\frac{2\phi_\varepsilon(0)\alpha_\varepsilon\mu_\varepsilon
}{1-\alpha_\varepsilon\mu_\varepsilon}\llVert
u_\varepsilon\rrVert _{L^1(\R
,\mrr x)},\qquad\varepsilon>0. \label{eq:bound2bis}
\end{eqnarray}
Here, $\mathfrak{N}_\varepsilon$ is defined as $\mathfrak{N}$ in
\eqref{eq:UBnonlin}, with $\phi_\varepsilon$, $u_\varepsilon$,
$\alpha_\varepsilon$ and $\mu_\varepsilon$
in place of $\phi$, $u$, $\alpha$ and $\mu$, respectively,
%
\begin{eqnarray*}
\delta^{(\varepsilon)}(u_\varepsilon)&:=&\int_{\R}u_\varepsilon
(t) \bigl(N_\varepsilon(\mrr t)-\lambda_\varepsilon(t) \,\mrr t
\bigr),
\\
\lambda_\varepsilon(t) &:=& \phi_\varepsilon \biggl(\int
_{(-\infty,t)}h_\varepsilon(t-s)N_\varepsilon(\mrr s)
\biggr),
\\
\delta_a^{(\varepsilon)}(u_\varepsilon)&:=&\int
_{\R}u_\varepsilon (t) \bigl(N_\varepsilon(\mrr t)-\widehat{\lambda_\varepsilon} \,\mrr t\bigr)
\end{eqnarray*}
and $\widehat{\lambda_\varepsilon}\in[\phi_\varepsilon(0),\phi
_\varepsilon(0)(1-\alpha_\varepsilon\mu_\varepsilon)^{-1}]$.

\item[(ii)] If, as $\varepsilon\to0$,
\begin{eqnarray}
\label{eq:lim1bis} \alpha_\varepsilon\mu_\varepsilon&\to&0,
\\
\label{eq:lim2bis} \phi_\varepsilon(0)\llVert u_\varepsilon\rrVert
_{L^2(\R,\mrr x)}^2&\to&1,
\\
\label{eq:lim3bis} \phi_\varepsilon(0)\llVert u_\varepsilon\rrVert
_{L^3(\R,\mrr x)}^3&\to&0,
\\
\label{eq:lim4bis} \sqrt{\phi_\varepsilon(0)}\alpha_\varepsilon
\mu_\varepsilon\bigl\llVert (u_\varepsilon)^2\bigr\rrVert
_{L^2(\R,\mrr x)}&\to&0,
\end{eqnarray}
then
\[
d_W\bigl(\delta^{(\varepsilon)}(u_\varepsilon),Z\bigr)\to0,
\qquad\mbox{as }\varepsilon\to0.
\]
If moreover, as $\varepsilon\to0$,
\begin{equation}
\label{eq:lim5bis} \phi_\varepsilon(0)\alpha_\varepsilon\mu_\varepsilon
\llVert u_\varepsilon\rrVert _{L^1(\R,\mrr x)}\to0,
\end{equation}
then
\[
d_W\bigl(\delta_a^{(\varepsilon)}(u_\varepsilon),Z
\bigr)\to0,\qquad\mbox{as }\varepsilon\to0.
\]
\end{longlist}
\end{Corollary}

We conclude this subsection with an example.

\begin{Example}\label{ex:1}
Let $I_\varepsilon$, $\varepsilon>0$, be a given family of bounded
Borel sets, $I_\varepsilon$ with Lebesgue measure $\ell_\varepsilon$,
and $\phi_\varepsilon:[0,\infty)\to[0,\infty)$, $\phi_\varepsilon
(0)>0$, $\varepsilon>0$,
a family of nondecreasing and Lipschitz continuous functions with
Lipschitz constant $\alpha_\varepsilon$.
Let $\mu_\varepsilon$, $\varepsilon>0$, be a collection of positive
numbers such that $\alpha_\varepsilon\mu_\varepsilon\subset(0,1)$,
$\varepsilon>0$,
and define the functions $h_\varepsilon(t):=\mu_\varepsilon
f_\varepsilon(t)$, $\varepsilon>0$, $t>0$, where
$f_\varepsilon$ is a locally bounded
probability density $($with respect to the Lebesgue measure$)$
on $(0,\infty)$. Hereafter, we consider
the family $N_\varepsilon$, $\varepsilon>0$, of stationary nonlinear
Hawkes processes with parameters $(\phi_\varepsilon,h_\varepsilon)$,
$\varepsilon>0$,
and the functions
\[
u_\varepsilon(t):=\frac{1}{\sqrt{\frac{\phi_\varepsilon(0)\ell
_\varepsilon}{1-\alpha_\varepsilon\mu_\varepsilon}}}\ind _{I_\varepsilon}(t),\qquad
\varepsilon>0, t\in\R.
\]
We have
\begin{eqnarray*}
\llVert u_\varepsilon\rrVert _{L^2(\R,\mrr x)}^2&=&
\frac{1-\alpha
_\varepsilon\mu_\varepsilon}{\phi_\varepsilon(0)}, \qquad\llVert u_\varepsilon\rrVert _{L^3(\R,\mrr x)}^3=
\biggl(\frac
{1-\alpha_\varepsilon\mu_\varepsilon}{\phi_\varepsilon(0)\ell
_\varepsilon} \biggr)^{3/2}\ell_\varepsilon,
\\
\bigl\llVert (u_\varepsilon)^2\bigr\rrVert _{L^2(\R,\mrr x)}&=&
\frac{1-\alpha
_\varepsilon\mu_\varepsilon}{\phi_\varepsilon(0)\sqrt{\ell
_\varepsilon}} \quad\mbox{and} \quad\llVert u_\varepsilon\rrVert
_{L^1(\R,\mrr x)}= \biggl(\frac{1-\alpha
_\varepsilon\mu_\varepsilon}{\phi_\varepsilon(0)\ell_\varepsilon
} \biggr)^{1/2}
\ell_\varepsilon.
\end{eqnarray*}
So by Corollary~\ref{cor:QCLTnonlinear}(i) we deduce
\begin{eqnarray}
d_W\bigl(\delta^{(\varepsilon)}(u_\varepsilon),Z\bigr)\leq
\mathfrak {N}_\varepsilon,\qquad\varepsilon>0\label{eq:RHSbis}
\end{eqnarray}
and
\begin{eqnarray}
d_W\bigl(\delta_a^{(\varepsilon)}(u_\varepsilon),Z
\bigr)\leq\mathfrak {N}_\varepsilon+ \frac{2\sqrt{\phi_\varepsilon(0)\ell_\varepsilon}\alpha
_\varepsilon\mu_\varepsilon}{\sqrt{1-\alpha_\varepsilon\mu
_\varepsilon}},\qquad\varepsilon>0,
\label{eq:RHS1bis}
\end{eqnarray}
where
\begin{eqnarray}
\mathfrak{N}_\varepsilon&:=&\sqrt{2/\pi}\alpha_\varepsilon\mu
_\varepsilon+\frac{2\sqrt{2/\pi}\alpha_\varepsilon\mu
_\varepsilon(2-\alpha_\varepsilon\mu_\varepsilon)}{
1-\alpha_\varepsilon\mu_\varepsilon}+\sqrt{\frac{1-\alpha
_\varepsilon\mu_\varepsilon}{\phi_\varepsilon(0)\ell_\varepsilon
}}
\nonumber
\\
&&{}+\frac{\alpha_\varepsilon\mu_\varepsilon}{\sqrt{\phi
_\varepsilon(0)\ell_\varepsilon(1-\alpha_\varepsilon\mu
_\varepsilon)}}.
\nonumber
\end{eqnarray}
If
\begin{equation}
\label{eq:Hyp1} \lim_{\varepsilon\to0}\phi_\varepsilon(0)
\ell_\varepsilon =+\infty\quad\mbox{and}\quad \lim_{\varepsilon\to0}
\alpha_\varepsilon\mu_\varepsilon=0,
\end{equation}
then one may easily check conditions \eqref{eq:lim1bis}, \eqref
{eq:lim2bis}, \eqref{eq:lim3bis} and \eqref{eq:lim4bis} and so
$d_W(\delta^{(\varepsilon)}(u_\varepsilon),Z)\to0$, as $\varepsilon
\to0$. To guarantee condition \eqref{eq:lim5bis} and, therefore,
$d_W(\delta_a^{(\varepsilon)}(u_\varepsilon),Z)\to0$,
as $\varepsilon\to0$, we need to suppose
\begin{equation}
\label{eq:Hyp2} \alpha_\varepsilon\mu_\varepsilon=o \biggl(
\frac{1}{\sqrt{\phi
_\varepsilon(0)\ell_\varepsilon}} \biggr),\qquad\mbox{as }\varepsilon\to0.
\end{equation}
Clearly, for specific choices of the sets $I_\varepsilon$,
$\varepsilon>0$, and the quantities $\phi_\varepsilon(0)$, $\alpha
_\varepsilon$ and
$\mu_\varepsilon$, we can provide the rate of convergence to zero of
the Wasserstein distances. For instance if, for $\varepsilon\in(0,1)$,
we take $I_\varepsilon=(0,1/\varepsilon)$, $\phi_\varepsilon(0)=\nu
$, $\alpha_\varepsilon=\alpha$, being $\nu$ and $\alpha$ positive
constants, and
$\mu_\varepsilon=\varepsilon$, then a straightforward computation
shows that, as $\varepsilon\to0$,
the right-hand sides of \eqref{eq:RHSbis} and \eqref{eq:RHS1bis}
converge to a positive constant when divided by $\varepsilon^{1/2}$,
and so
\[
d_W\bigl(\delta^{(\varepsilon)}(u_\varepsilon),Z
\bigr)=d_W\bigl(\delta _a^{(\varepsilon)}(u_\varepsilon),Z
\bigr) =O(\sqrt{\varepsilon}),\qquad\mbox{as }\varepsilon\to0.
\]
The bounds \eqref{eq:RHSbis} and \eqref{eq:RHS1bis} may be used to
construct confidence intervals for $\delta^{(\varepsilon
)}(u_\varepsilon)$
and $\delta_a^{(\varepsilon)}(u_\varepsilon)$, respectively. For
instance, let $F_X$ denote the distribution function of a random
variable $X$,
assume \eqref{eq:Hyp1}, choose $b_i^{(\beta)}$, $i=1,2$, $\beta\in
(0,1/2)$, so that $b_1^{(\beta)}< b_2^{(\beta)}$ and
\begin{equation}
\label{eq:confid1} F_Z\bigl(b_1^{(\beta)}\bigr)\leq
\beta/2\quad\mbox{and}\quad F_Z\bigl(b_2^{(\beta
)}
\bigr)\geq1-\frac{\beta}{2}
\end{equation}
and choose $\varepsilon>0$ so small that $2\sqrt{\mathfrak
{N}_\varepsilon}\leq\beta/2$, then
\begin{equation}
\label{eq:confid2} P\bigl(b_1^{(\beta)}<\delta^{(\varepsilon)}(u_\varepsilon)
\leq b_2^{(\beta)}\bigr)\geq1-2\beta.
\end{equation}
Indeed, by \eqref{eq:RHSbis} and the inequality
\[
\sup_{x\in\R}\bigl\llvert F_X(x)-F_Z(x)
\bigr\rrvert \leq2\sqrt{d_W(X,Z)}
\]
(see, e.g., \cite{chen}), we have
\[
\bigl\llvert F_{\delta^{(\varepsilon)}(u_\varepsilon)}\bigl(b_i^{(\beta
)}
\bigr)-F_Z\bigl(b_i^{(\beta)}\bigr)\bigr\rrvert
\leq2\sqrt{\mathfrak{N}_\varepsilon}\leq \beta/2,\qquad i=1,2.
\]
Relation \eqref{eq:confid2} easily follows by this latter inequality
and \eqref{eq:confid1}.
\end{Example}

\subsection{Proofs of Theorems \texorpdfstring{\protect\ref{thm:nonlinGauss}}{4.1} and \texorpdfstring{\protect\ref{thm:nonlinGaussappr}}{4.3}}\label{sec:proofs1}

\mbox{}

\begin{pf*}{Proof of Theorem~\ref{thm:nonlinGauss}}
We divide the proof in two main steps. In the first step, we prove
\begin{eqnarray}\label{eq:thm41}
d_W\bigl(\delta(u),Z\bigr)&\leq&\sqrt{2/\pi}\mathrm{E} \biggl[\biggl
\llvert 1-\int_{\R}\bigl\llvert u(t)\bigr\rrvert
^2\lambda(t) \,\mrr t\biggr\rrvert \biggr]\nonumber
\\
&&{}+\lambda\llVert u\rrVert
_{L^3(\R,\mrr x)}^3
\nonumber\\[-8pt]\\[-8pt]\nonumber
&&{}+\frac{2\lambda\alpha\mu}{1-\alpha\mu}\sqrt{2/\pi}\llVert u\rrVert _{L^2(\R,\mrr x)}^{2}
\\
&&{}
+\frac{\lambda\alpha\mu}{1-\alpha\mu}\llVert u\rrVert _{L^2(\R,\mathrm
{d}x)}\bigl\llVert u^2
\bigr\rrVert _{L^2(\R,\mrr x)}.\nonumber 
\end{eqnarray}
In the second step, we complete the proof. If $u\notin L^2(\R,\mathrm
{d}x)\cap L^3(\R,\mrr x)\cap L^4(\R,\mrr x)$, then the
claim is clearly true.
So we shall assume $u\in L^2(\R,\mrr x)\cap L^3(\R,\mathrm
{d}x)\cap L^4(\R,\mrr x)$.

\begin{longlist}
\item[\textit{Step}  1:  Proof  of  \eqref{eq:thm41}.] Hereafter, for
ease of notation, we write $\varphi_t(\overline{N}\mid_{(-\infty,t)})$
in place of
$\varphi(t,\overline{N}\mid_{(-\infty,t)})$, $t\in\R$. By \eqref
{eq:functional}, for $s,t\in\R$, we have
\begin{eqnarray}
\lambda(s)&=&\varphi_s(\overline{N}\mid_{(-\infty,s)})
\nonumber
\\
&=&\phi \biggl(\int_{(-\infty,s)\times\R_+}
\ind_{u\leq
t}h(s-u)\ind_{(0,\varphi_u(\overline{N}\mid_{(-\infty,u)})]}(v) \overline{N}(\mrr u\times \mrr v)
\nonumber
\\
&&{} +\int_{(-\infty,s)\times\R_+}\ind_{u>t}h(s-u)\ind_{(0,\varphi
_u(\overline{N}\mid_{(-\infty,u)})]}(v)
\overline{N}(\mrr u\times \mrr v) \biggr).
\nonumber
\end{eqnarray}
We shall show later on that $h\geq0$ and $\phi$ nondecreasing imply
\begin{equation}
\label{eq:Dpos} D_{(t,z)}\lambda(s)\geq0,\qquad\mbox{for }s,t\in\R, s>t
\mbox{ and }z\in\R_+.
\end{equation}
So, by the Lipschitz continuity of $\phi$, for $s,t\in\R$, $s>t$,
and $z\in\R_+$, we have
\begin{eqnarray} \label{eq:Lip}
0&\leq& D_{(t,z)}\lambda(s) 
\nonumber
\\
&\leq&\alpha \biggl(\int_{(-\infty,s)\times\R_+}\ind_{u\leq
t}h(s-u)
\ind_{(0,\varphi_u((\overline{N}+\varepsilon
_{(t,z)})\mid_{(-\infty,u)})]}(v) (\overline{N}+\varepsilon_{(t,z)}) (\mrr u
\times\mathrm {d}v)
\nonumber
\\
&&{} +\int_{(-\infty,s)\times\R_+}\ind_{u>t}h(s-u)\ind_{(0,\varphi
_u((\overline{N}+\varepsilon_{(t,z)})\mid_{(-\infty,u)})]}(v)
(\overline{N}+\varepsilon_{(t,z)}) (\mrr u\times\mathrm {d}v)
\nonumber
\\
&&{} -\int_{(-\infty,s)\times\R_+}\ind_{u\leq t}h(s-u)\ind_{(0,\varphi
_u(\overline{N}\mid_{(-\infty,u)})]}(v)
\overline{N}(\mrr u\times \mrr v)
\nonumber
\\
&&{} -\int_{(-\infty,s)\times\R_+}\ind_{u>t}h(s-u)\ind_{(0,\varphi
_u(\overline{N}\mid_{(-\infty,u)})]}(v)
\overline{N}(\mrr u\times \mrr v) \biggr)
\nonumber
\\
&=&\alpha \biggl(\int_{(-\infty,s)\times\R_+}\ind_{u\leq
t}h(s-u)
\ind_{(0,\varphi_u(\overline{N}\mid_{(-\infty,u)})]}(v) (\overline{N}+\varepsilon_{(t,z)}) (\mrr u
\times\mathrm {d}v)
\nonumber
\\
&&{} +\int_{(-\infty,s)\times\R_+}\ind_{u>t}h(s-u)\ind_{(0,\varphi
_u(\overline{N}\mid_{(-\infty,u)}+\varepsilon_{(t,z)})]}(v)
(\overline{N}+\varepsilon_{(t,z)}) (\mrr u\times\mathrm {d}v)
\nonumber
\\
&&{} -\int_{(-\infty,s)\times\R_+}\ind_{u\leq t}h(s-u)\ind_{(0,\varphi
_u(\overline{N}\mid_{(-\infty,u)})]}(v)
\overline{N}(\mrr u\times \mrr v)
\nonumber\\[-8pt]\\[-8pt]\nonumber
&&{} -\int_{(-\infty,s)\times\R_+}\ind_{u>t}h(s-u)\ind_{(0,\varphi
_u(\overline{N}\mid_{(-\infty,u)})]}(v)
\overline{N}(\mrr u\times \mrr v) \biggr)
\nonumber
\\
&=&\alpha \biggl(h(s-t)\ind_{(0,\varphi_t(\overline{N}\mid_{(-\infty,t)})]}(z)\nonumber
\\
&&{} + \int_{(-\infty,s)\times\R_+}\ind_{u>t}h(s-u)\ind_{(0,\varphi
_u(\overline{N}\mid_{(-\infty,u)}+\varepsilon_{(t,z)})]}(v)
\overline{N}(\mrr u\times \mrr v)
\nonumber
\\
&&{} -\int_{(-\infty,s)\times\R_+}\ind_{u>t}h(s-u)\ind_{(0,\varphi
_u(\overline{N}\mid_{(-\infty,u)})]}(v)
\overline{N}(\mrr u\times \mrr v) \biggr)
\nonumber
\\
&=&\alpha \biggl(h(s-t)\ind_{(0,\varphi_t(\overline{N}\mid_{(-\infty
,t)})]}(z)
\nonumber
\\
&&{} +\int_{(t,s)\times\R_+}h(s-u) \bigl(\ind_{(0,\varphi_u(\overline
{N}\mid_{(-\infty,u)}+\varepsilon_{(t,z)})]}(v) \nonumber
\\
&&{} -
\ind_{(0,\varphi_u(\overline{N}\mid_{(-\infty,u)})]}(v)\bigr) \overline{N}(\mrr u\times \mrr v) \biggr)
\nonumber
\\
&=&\alpha \biggl(h(s-t)\ind_{(0,\varphi_t(\overline{N}\mid_{(-\infty,t)})]}(z) +\int_{(t,s)}h(s-u)N_{(t,z)}(\mrr u) \biggr),\nonumber
\end{eqnarray}
where $N_{(t,z)}$ is the point process on $(t,\infty)$ defined by
\[
N_{(t,z)}(\mrr u):=\overline{N}\bigl(\mrr u\times\bigl(\varphi
_u(\overline{N}\mid_{(-\infty,u)}), \varphi_u(
\overline{N}\mid_{(-\infty,u)}+\varepsilon_{(t,z)})\bigr]\bigr).
\]
The processes $\{\varphi_u(\overline{N}\mid_{(-\infty,u)}+\varepsilon
_{(t,z)})\}_{u>t}$ and
$\{\varphi_u(\overline{N}\mid_{(-\infty,u)})\}_{u>t}$ are $\{\mathcal
{F}_u^{\overline{N}}\}_{u>t}$-predictable, and
we shall check later on that the mapping
\begin{equation}
\label{eq:locbound} (t,\infty)\ni u\mapsto\mathrm{E}\bigl[D_{(t,z)}\lambda(u)
\bigr]\in\R\mbox{ is locally bounded.}
\end{equation}
%
Therefore,
\[
\int_a^b D_{(t,z)}\lambda(u)
\,\mrr u<\infty\qquad\mbox{a.s., for any }a,b>t.
\]
Consequently, by Lemma~\ref{le:PoissEmb} we have that $N_{(t,z)}$ has
$\{\mathcal{F}_u^{\overline{N}}\}_{u>t}$-stochastic intensity $\{
D_{(t,z)}\lambda(u)\}_{u>t}$.
Taking the mean in \eqref{eq:Lip}, we deduce
\begin{eqnarray}
\mathrm{E}\bigl[D_{(t,z)}\lambda(s)\bigr]&\leq&\alpha \biggl(h(s-t)P\bigl(
\lambda (t)\geq z\bigr) +\int_{t}^{s}h(s-u)
\mathrm{E}\bigl[D_{(t,z)}\lambda(u)\bigr] \,\mathrm {d}u \biggr).
\nonumber
\end{eqnarray}
Extending the definition of $h$ for nonpositive times as $h(t)=0$,
$t\leq0$, we rewrite the above inequality as
\begin{eqnarray}
q_{(t,z)}(s)\leq p_{(t,z)}(s)+r*q_{(t,z)}(s),\qquad s,t\in
\R, z\in\R _+,
\nonumber
\end{eqnarray}
where for ease of notation we set $q_{(t,z)}(s):=\mathrm
{E}[D_{(t,z)}\lambda(s)]$, $p_{(t,z)}(s):=\alpha h(s-t)P(\lambda
(t)\geq z)$,
$r(s):=\alpha h(s)$ and $*$ denotes the convolution product between
functions. Iterating this inequality, we deduce, for $n\geq1$,
\begin{eqnarray}
q_{(t,z)}(s)\leq\sum_{i=0}^{n-1}p_{(t,z)}*r^{*i}(s)+q_{(t,z)}*r^{*n}(s),
\qquad s,t\in\R, z\in\R_+,
\nonumber
\end{eqnarray}
where $r^{*0}$ is by definition the Dirac delta function. By \eqref
{eq:locbound} and the stability condition $\alpha\mu<1$,
we deduce $q_{(t,z)}*r^{*n}(s)\to0$, as $n\to\infty$, for any
$t,s\in\R$, $z\in\R_+$. Indeed, for some constant $C_{t,z,s}>0$,
\begin{eqnarray}
q_{(t,z)}*r^{*n}(s)&=&\int_{\R}r^{*n}(s-u)q_{(t,z)}(u)\,
\mathrm {d}u =\int_{t}^{s}r^{*n}(s-u)q_{(t,z)}(u)
\,\mrr u
\nonumber
\\
&\leq& C_{t,z,s}\int_{\R}r^{*n}(s-u)
\,\mrr u
\nonumber
\\
&\leq& C_{t,z,s}(\alpha\mu)^{n},
\nonumber
\end{eqnarray}
where the latter inequality follows by a standard property of
convolutions; see, for example, Theorem IV.15 in \cite{brezis}.
Therefore,
\begin{eqnarray}\label{eq:qdis}
q_{(t,z)}(s)&\leq&\sum_{i\geq0}p_{(t,z)}*r^{*i}(s)\nonumber
\\
&=& P \bigl(\lambda(t)\geq z\bigr)\sum_{i\geq0}
\alpha^{i+1}\int_{\R}h(s-u-t) h^{*i}(u)\,
\mathrm {d}u
\\
&=&P\bigl(\lambda(t)\geq z\bigr)\sum_{i\geq1}
\alpha^{i}h^{*i}(s-t),\qquad s,t\in\R, z\in
\R_+.\nonumber
\end{eqnarray}
Consequently, for any $f,g$ integrable and square integrable, defining
$\check{f}(x):=f(-x)$, by the Cauchy--Schwarz inequality and the
properties of the
convolution product (see again Theorem IV.15 in \cite{brezis}), we have
\begin{eqnarray}
&&\int_{\R\times\R_+}\bigl\llvert f(t)\bigr\rrvert \biggl(\int
_{t}^{+\infty} \bigl\llvert g(s)\bigr\rrvert \mathrm{E}
\bigl[\ind_{(0,\lambda(t)]}(z)D_{(t,z)}\lambda(s)\bigr] \,\mrr s \biggr)
\,\mrr t\,\mrr z
\nonumber
\\
&&\qquad \leq 
\lambda\sum_{i\geq1}\alpha^i\int
_{\R}\bigl\llvert f(t)\bigr\rrvert \bigl(\check
{h^{*i}}*\llvert g\rrvert \bigr) (t) \,\mrr t\label{eq:intint}
\\
&&\qquad \leq\lambda\llVert f\rrVert _{L^2(\R,\mrr x)}\sum_{i\geq1}
\alpha^i\bigl\llVert \check{h^{*i}}*\llvert g\rrvert \bigr
\rrVert _{L^2(\R,\mrr x)}
\nonumber
\\
&&\qquad \leq\lambda\llVert f\rrVert _{L^2(\R,\mrr x)}\llVert g\rrVert _{L^2(\R,\mrr x)}
\sum_{i\geq1}\alpha^i\bigl\llVert
\check{h^{*i}}\bigr\rrVert _{L^1(\R,\mathrm
{d}x)}\nonumber
\\
&&\qquad \leq\lambda\llVert f\rrVert _{L^2(\R,\mrr x)}\llVert g\rrVert _{L^2(\R,\mrr x)}
\sum_{i\geq1}\alpha^i\mu^i
\nonumber\\[-8pt]\label{eq:intintint} \\[-8pt]\nonumber
&&\qquad =
\llVert f\rrVert _{L^2(\R,\mrr x)}\llVert g\rrVert _{L^2(\R,\mrr x)}
\frac{\lambda\alpha\mu}{1-\alpha\mu}.
\end{eqnarray}
Assume for the moment that we may apply Theorem~\ref{thm:normalskor},
then by \eqref{eq:boundgauss} and the above inequality (applied first
with $f=g=u$ and then
with $f=u$ and $g=u^2$) we easily deduce \eqref{eq:thm41}.

It remains to prove \eqref{eq:Dpos}, \eqref{eq:locbound} and to check
the assumptions of Theorem~\ref{thm:normalskor}.

We first prove \eqref{eq:locbound}. Let $(t,z)\in\R\times\R_+$ be fixed.
Since
\[
\mathrm{E}\bigl[\lambda(u)\bigr]=\mathrm{E}\bigl[\varphi_u(\overline
{N}\mid_{(-\infty,u)})\bigr]=\lambda\in(0,\infty)
\]
for any $u\in\R$, to show
\eqref{eq:locbound} it suffices to prove that
the map
\[
u\mapsto\mathrm{E}\bigl[\varphi_u(\overline{N}\mid_{(-\infty
,u)}+\varepsilon_{(t,z)})\bigr]
\]
is locally bounded on $(t,\infty)$.
We define recursively the processes $\lambda_{(t,z)}^{\prime(0)}\equiv0$,
\begin{eqnarray}\label{eq:lambdanprimo}
N_{(t,z)}^{\prime(n)}(\mrr s)&=&\overline{N}\bigl(\mrr s\times
\bigl(0,\lambda_{(t,z)}^{\prime(n)}(s)\bigr]\bigr),
\nonumber
\\
\lambda_{(t,z)}^{\prime(n+1)}(s)&=&\phi \biggl(\int_{(-\infty
,s)}h(s-u)\bigl(N_{(t,z)}^{\prime(n)}(\mrr u)
\\
&&{} +\varepsilon_{(t,z)}
\bigl(\mathrm {d}u\times\bigl(0,\lambda_{(t,z)}^{\prime(n)}(u)\bigr]\bigr)
\bigr) \biggr),\qquad n\geq0,  s>t.\nonumber
\end{eqnarray}
We are going to check by induction that, for any $n\geq0$,
\begin{equation}
\label{eq:ind1} \int_a^b\lambda_{(t,z)}^{\prime(n)}(s)
\,\mrr s<\infty,\qquad\mbox {a.s., for any }a,b>t
\end{equation}
and $\{\lambda_{(t,z)}^{\prime(n)}(s)\}_{s>t}$ is $\{\mathcal
{F}_s^{\overline{N}}\}_{s>t}$-predictable. The basis of the induction
is clearly verified.
So assume the claim for $\lambda_{(t,z)}^{\prime(n)}$ and let $\{
T_{(t,z),m}^{\prime(n)}\}_{m\in\mathbb{Z}}$ be the points of
$N_{(t,z)}^{\prime(n)}$ on $(t,\infty)$.
By Lemma~\ref{le:PoissEmb}, we have that $N_{(t,z)}^{\prime(n)}$ has $\{
\mathcal{F}_s^{\overline{N}}\}_{s>t}$-stochastic intensity $\{\lambda
_{(t,z)}^{\prime(n)}(s)\}_{s>t}$.
By the Lipschitz property of $\phi$ and the nonnegativity of $h$,
we deduce
\begin{eqnarray}\label{eq:liplambda'}
\lambda_{(t,z)}^{\prime(n+1)}(s)&\leq&\phi(0)+\alpha\int
_{(-\infty
,s)}h(s-u)N_{(t,z)}^{\prime(n)}(\mrr u)\nonumber
\\
&&{} +\alpha
\int_{(-\infty
,s)}h(s-u)\varepsilon_{(t,z)}\bigl(\mrr u\times
\bigl(0,\lambda _{(t,z)}^{\prime(n)}(u)\bigr]\bigr)
\nonumber
\\
&=&\phi(0)+\alpha\sum_{m\in\mathbb{Z}}h\bigl(s-T_{(t,z),m}^{\prime(n)}
\bigr)\ind _{(-\infty,s)}\bigl(T_{(t,z),m}^{\prime(n)}\bigr)
\nonumber\\[-8pt]\\[-8pt]\nonumber
&&{} +\alpha\int_{(-\infty,s)\times\R_+}h(s-u)\ind_{(0,\lambda
_{(t,z)}^{\prime(n)}(u)]}(v)
\varepsilon_{(t,z)}(\mrr u\times\mathrm {d}v)
\nonumber
\\
&=&\phi(0)+\alpha\sum_{m\in\mathbb{Z}}h\bigl(s-T_{(t,z),m}^{\prime(n)}
\bigr)\ind _{(-\infty,s)}\bigl(T_{(t,z),m}^{\prime(n)}\bigr) \nonumber
\\
&&{} +\alpha
h(s-t)\ind_{(0,\lambda_{(t,z)}^{\prime(n)}(t)]}(z).\nonumber
\end{eqnarray}
Integrating over the finite interval $(a,b)\subset(t,\infty)$, we have
\begin{eqnarray}
\int_a^b\lambda_{(t,z)}^{\prime(n+1)}(s)
\,\mrr s &\leq&\phi(0) (b-a)+\alpha\sum_{m\in\mathbb{Z}}\ind
\bigl\{ t<T_{(t,z),m}^{\prime(n)}<b\bigr\}\int_{0\vee
(a-T_{(t,z),m}^{\prime(n)})}^{b-T_{(t,z),m}^{\prime(n)}}h(u)
\,\mrr u
\nonumber
\\
&&{}+\alpha\int_a^b h(s-t) \,\mrr s,
\nonumber
\end{eqnarray}
and this latter quantity is finite since $h$ is integrable and
$N_{(t,z)}^{\prime(n)}$ has an a.s. finite number of points\vspace*{1pt} in any bounded
interval of $(t,\infty)$
[due to \eqref{eq:ind1}]. Moreover, the process $\{\lambda
_{(t,z)}^{\prime(n+1)}(s)\}_{s>t}$ is $\{\mathcal{F}_s^{\overline N}\}
_{s>t}$-predictable. Indeed,
\[
\lambda_{(t,z)}^{\prime(n+1)}(s)=\phi \biggl(h(s-t)\ind_{(0,\lambda
_{(t,z)}^{\prime(n)}(t)]}(z)+
\int_{(-\infty
,s)}h(s-u)N_{(t,z)}^{\prime(n)}(\mrr u)
\biggr)
\]
and the processes
\begin{equation}
\label{eq:predictability} \bigl\{h(s-t)\ind_{(0,\lambda_{(t,z)}^{\prime(n)}(t)]}(z)\bigr\}_{s>t},\qquad
\biggl\{\int_{(-\infty,s)}h(s-u)N_{(t,z)}^{\prime(n)}(\mrr u) \biggr\}_{s>t}
\end{equation}
are $\{\mathcal{F}_s^{\overline{N}}\}_{s>t}$-predictable. To justify
the predictability of the first process in \eqref{eq:predictability},
one may first note that it is $\{\mathcal{F}_s^{\overline{N}}\}
_{s>t}$-adapted [since $\lambda_{(t,z)}^{\prime(n)}(t)$ is $\mathcal
{F}_t^{\overline{N}}$-measurable
and $h$ is deterministic] and then conclude by applying, for example,
Theorem T34 in~\cite{bremaud}.
To justify the predictability of the second process in \eqref
{eq:predictability}, one
notes that
it is left-continuous and $\{\mathcal{F}_s^{\overline{N}}\}
_{s>t}$-adapted. The induction is therefore completed and
by Lemma~\ref{le:PoissEmb}, for any $n\geq0$ and $(t,z)\in\R\times
\R_+$, the point process $N_{(t,z)}^{\prime(n)}$ on $(t,\infty)$
has $\{\mathcal{F}_s^{\overline{N}}\}_{s>t}$-stochastic intensity $\{
\lambda_{(t,z)}^{\prime(n)}(s)\}_{s>t}$. For fixed $(t,z)\in\R\times\R_+$,
since $h$ is nonnegative and $\phi$ is nondecreasing, we have that
$\lambda_{(t,z)}^{\prime(n)}(s,\omega)$ and $N_{(t,z)}^{\prime(n)}(C)(\omega)$
increase with $n\geq0$,
for all $\omega$, $s>t$ and Borel sets $C\subseteq(t,\infty)$. So
the limiting processes $\{\lambda_{(t,z)}^{\prime(\infty)}(s)\}_{s>t}$
and $N_{(t,z)}^{\prime(\infty)}$ are\vspace*{1pt} defined for all $\omega$. Setting
$h\equiv0$ on $(-\infty,0]$, by \eqref{eq:liplambda'}, for any
$n\geq0$, $s,t\in\R$ and $z\in\R_+$, we have
\[
\lambda_{(t,z)}^{\prime(n+1)}(s)\leq\phi(0)+\alpha h(s-t)+\alpha\int
_{\R}h(s-u)N_{(t,z)}^{\prime(n)}(\mrr u).
\]
Taking the mean over this inequality, we have
\[
q_{(t,z)}^{\prime(n+1)}(s)\leq p_t'(s)+r*q_{(t,z)}^{\prime(n)}(s),
\]
where for ease of notation we set $q_{(t,z)}^{\prime(n)}(s):=\mathrm
{E}[\lambda_{(t,z)}^{\prime(n)}(s)]$, $p_t'(s):=\phi(0)+r(s-t)$ and the
function $r$ is defined as above.
Iterating this latter inequality and using that $q_{(t,z)}^{\prime(0)}\equiv
0$, we deduce
\begin{eqnarray}
q_{(t,z)}^{\prime(n+1)}(s)&\leq&\sum_{i\geq0}r^{*i}*p_t'(s)=
\phi(0)\sum_{i\geq0}\bigl\llVert r^{*i}\bigr
\rrVert _{L^1(\R,\mrr x)}+\sum_{i\geq
1}r^{*i}(s-t).
\nonumber
\end{eqnarray}
Passing to the limit as $n\to\infty$, by the monotone convergence
theorem, a standard property of the convolution and the stability
condition $\alpha\mu<1$, we have
\begin{eqnarray}
\mathrm{E}\bigl[\lambda_{(t,z)}^{\prime(\infty)}(s)\bigr]&=:&q_{(t,z)}^{\prime(\infty)}(s)\leq\phi(0)\sum_{i\geq0}\bigl\llVert r^{*i}
\bigr\rrVert _{L^1(\R,\mrr x)}+\sum_{i\geq1}r^{*i}(s-t)
\nonumber
\\
&\leq&\frac{\phi(0)}{1-\alpha\mu}+\sum_{i\geq
1}r^{*i}(s-t)
\nonumber
\\
&=&\frac{\phi(0)}{1-\alpha\mu}+\sum_{i\geq1}\alpha^i
\int_{\R
}h^{* i-1}(s-t-u)h(u) \,\mrr u
\nonumber
\\
&\leq&\frac{\phi(0)}{1-\alpha\mu}+\sum_{i\geq1}
\alpha^i\int_{0}^{s-t}h^{* i-1}(s-t-u)h(u)
\,\mrr u
\nonumber
\\
&\leq&\frac{\phi(0)}{1-\alpha\mu}+\alpha(1-\alpha\mu)^{-1}
\ind _{s>t}\max_{u\in[0,s-t]}h(u),
\nonumber
\end{eqnarray}
and so by the local boundedness of $h$ we have
\begin{equation}
\label{eq:locbouninfty} \max_{s\in[a,b]}q_{(t,z)}^{\prime(\infty)}(s)<
\infty,\qquad\mbox{for any }a<b, a,b>t.
\end{equation}
%
In particular,
\[
\int_a^b\lambda_{(t,z)}^{\prime(\infty)}(s)
\,\mrr s<\infty\qquad \mbox{a.s., for any }a<b, a,b>t.
\]
Moreover, $\{\lambda_{(t,z)}^{\prime(\infty)}(s)\}_{s>t}$ is $\{\mathcal
{F}_s^{\overline{N}}\}_{s>t}$-predictable as limit of
$\{\mathcal{F}_s^{\overline{N}}\}_{s>t}$-predictable processes and
\begin{eqnarray}
N_{(t,z)}^{\prime(\infty)}(\mrr s)&=&\overline{N}\bigl(\mrr s\times
\bigl(0,\lambda_{(t,z)}^{\prime(\infty)}(s)\bigr]\bigr),\qquad s>t.
\nonumber
\end{eqnarray}
So by Lemma~\ref{le:PoissEmb} $N_{(t,z)}^{\prime(\infty)}$ has $\{\mathcal
{F}_s^{\overline{N}}\}_{s>t}$-stochastic intensity $\{\lambda
_{(t,z)}^{\prime(\infty)}(s)\}_{s>t}$.
Taking the limit as $n\to\infty$ in \eqref{eq:lambdanprimo}, we have
\begin{eqnarray}\label{eq:lambdainfeq}
\lambda_{(t,z)}^{\prime(\infty)}(s) &=& \phi \biggl(\int
_{(-\infty
,s)}h(s-u) \bigl(N_{(t,z)}^{\prime(\infty)}(\mrr u)
\nonumber\\[-8pt]\\[-8pt]\nonumber
&&{} +\varepsilon _{(t,z)}\bigl(\mrr u\times\bigl(0,
\lambda_{(t,z)}^{\prime(\infty
)}(u)\bigr]\bigr)\bigr) \biggr), \qquad s>t.
\end{eqnarray}
Therefore,
\[
\mathrm{E}\bigl[\varphi_s\bigl((\overline{N}+\varepsilon_{(t,z)})
\mid_{(-\infty,s)}\bigr)\bigr]= \mathrm{E}\bigl[\varphi_s(
\overline{N}\mid_{(-\infty,s)}+\varepsilon_{(t,z)})\bigr] =\mathrm{E}
\bigl[\lambda_{(t,z)}^{\prime(\infty)}(s)\bigr],\qquad s>t
\]
and \eqref{eq:locbound} follows by \eqref{eq:locbouninfty}.

We now prove \eqref{eq:Dpos}. Let $(t,z)\in\R\times\R_+$ be fixed.
We define recursively the processes $\lambda_{(t,z)}^{\prime\prime(0)}\equiv0$,
\begin{eqnarray}
N_{(t,z)}^{\prime\prime(n)}(\mrr s)&=&\overline{N}\bigl(\mrr s\times
\bigl(0,\lambda_{(t,z)}^{\prime\prime(n)}(s)\bigr]\bigr),
\nonumber
\\
\lambda_{(t,z)}^{\prime\prime(n+1)}(s)&=&\phi \biggl(\int_{(-\infty
,s)}h(s-u)N_{(t,z)}^{\prime\prime(n)}(\mrr u) \biggr),\qquad n\geq0, s>t
\nonumber
\end{eqnarray}
and note that since $h\geq0$ and $\phi$ is nondecreasing we have
\begin{equation}
\label{eq:ineqn} \lambda_{(t,z)}^{\prime\prime(n)}(s,\omega)\leq\lambda
_{(t,z)}^{\prime(n)}(s,\omega),\qquad\mbox{for all }\omega, n\geq0\mbox
{ and }s>t,
\end{equation}
where $\lambda_{(t,z)}^{\prime(0)}\equiv0$ and $\lambda
_{(t,z)}^{\prime(n+1)}$, $n\geq0$, is defined by \eqref{eq:lambdanprimo}.
Arguing as above,
we have that, for fixed $(t,z)\in\R\times\R_+$, $\lambda
_{(t,z)}^{\prime\prime(n)}(s,\omega)$ and $N_{(t,z)}^{\prime\prime(n)}(C)(\omega)$
increase with $n\geq0$,
for all $\omega$, $s>t$ and Borel sets $C\subseteq(t,\infty)$, and
the limiting processes $\{\lambda_{(t,z)}^{\prime\prime(\infty)}(s)\}_{s>t}$
and $N_{(t,z)}^{\prime\prime(\infty)}$ are such that
\[
\int_a^b\lambda_{(t,z)}^{\prime\prime(\infty)}(s)
\,\mrr s<\infty\qquad \mbox{a.s., for any }a<b, a,b>t
\]
$\{\lambda_{(t,z)}^{\prime\prime(\infty)}(s)\}_{s>t}$ is $\{\mathcal
{F}_s^{\overline{N}}\}_{s>t}$-predictable and
\begin{eqnarray*}
N_{(t,z)}^{\prime\prime(\infty)}(\mrr s)&=&\overline{N}\bigl(\mrr s\times
\bigl(0,\lambda_{(t,z)}^{\prime\prime(\infty)}(s)\bigr]\bigr),\qquad s>t,
\nonumber
\\
\lambda_{(t,z)}^{\prime\prime(\infty)}(s)&=&\phi \biggl(\int
_{(-\infty
,s)}h(s-u)N_{(t,z)}^{\prime\prime(\infty)}(\mrr u)
\biggr), \qquad s>t.
\end{eqnarray*}
Inequality \eqref{eq:Dpos} follows noticing that taking the limit
as $n\to\infty$ in \eqref{eq:ineqn} we have
\[
\lambda_{(t,z)}^{\prime\prime(\infty)}(s,\omega)\leq\lambda _{(t,z)}^{\prime(\infty)}(s,
\omega),\qquad\mbox{for almost all }\omega \mbox{ and any }s>t.
\]

We now check the assumptions of Theorem~\ref{thm:normalskor}. Since
$N$ is stationary with a finite intensity and $u$ is integrable and
square integrable,
conditions \eqref{eq:L1} and \eqref{eq:L2} hold. Arguing similarly to
\eqref{eq:intint}, for an integrable function $g$ we have
\begin{eqnarray}
\int_{\R\times\R_+} \biggl(\int_t^\infty
\bigl\llvert g(s)\bigr\rrvert \mathrm {E}\bigl[D_{(t,z)}\lambda(s)\bigr]
\,\mrr s \biggr) \,\mrr t\,\mrr z 
&
\leq&\lambda\sum_{i\geq1}\alpha^i\int
_{\R}\check {h^{*i}}*\llvert g\rrvert (t)
\,\mrr t
\nonumber
\\
&\leq&\frac{\lambda\alpha\mu}{1-\alpha\mu}\llVert g\rrVert _{L^1(\R,\mathrm
{d}x)}.
\nonumber
\end{eqnarray}
Taking first $g=u$ and then $g=u^2$, one then has that conditions
\eqref{eq:L3} and \eqref{eq:L4} are satisfied. Condition \eqref{eq:L5}
follows by \eqref{eq:intintint} with $f=u^2$ and $g=u$. Finally, the
local integrability on $(t,\infty)$ of
the random function $D_{(t,z)}\lambda(\cdot)$ is a consequence of
\eqref{eq:locbound}.
The proof is complete.

\item[\textit{Step}  2:  Proof  of  \eqref{eq:thm41explicit}.] We have
\begin{eqnarray}\label{eq:secondmom}
&& \mathrm{E} \biggl[\biggl\llvert 1-\int_{\R}\bigl\llvert u(t)
\bigr\rrvert ^2\lambda(t)\, \mathrm {d}t\biggr\rrvert \biggr]
\nonumber\\[-8pt]\\[-8pt]\nonumber
&&\qquad \leq\bigl
\llvert 1-\lambda\llVert u\rrVert _{L^2(\R,\mrr x)}^2\bigr\rrvert +\int
_{\R}\bigl\llvert u(t)\bigr\rrvert ^2\mathrm{E}
\bigl[\bigl\llvert \lambda(t)-\lambda\bigr\rrvert \bigr] \,\mrr t.
\end{eqnarray}
By \eqref{eq:disstocint} and \eqref{eq:disintens}, it follows
\begin{eqnarray}
\bigl\llvert \lambda(t)-\lambda\bigr\rrvert &\leq& 
\max \biggl\{\alpha\int_{(-\infty,t)}h(t-s)N(\mrr s),
\frac
{\phi(0)\alpha\mu}{1-\alpha\mu} \biggr\}\qquad\mbox{a.s., for all }t\in\R.
\nonumber
\end{eqnarray}
Taking the expectation on this relation and using the rightmost
inequality in \eqref{eq:disintens}, we have
\begin{equation}
\label{eq:L1lambda} \mathrm{E}\bigl[\bigl\llvert \lambda(t)-\lambda\bigr\rrvert \bigr]
\leq\frac{2\phi(0)\alpha\mu
}{1-\alpha\mu}\qquad\mbox{a.s., for all }t\in\R.
\end{equation}
Combining this latter inequality with \eqref{eq:secondmom} and \eqref
{eq:disintens}, we deduce
\begin{eqnarray} \label{eq:secondmom2}
&& \mathrm{E} \biggl[\biggl\llvert 1-\int_{\R}\bigl\llvert u(t)
\bigr\rrvert ^2\lambda(t) \,\mathrm {d}t\biggr\rrvert \biggr]\nonumber
\\
&&\qquad \leq \bigl
\llvert 1-\lambda\llVert u\rrVert _{L^2(\R,\mrr x)}^2\bigr\rrvert +
\frac{2\phi(0)\alpha\mu}{1-\alpha\mu}\llVert u\rrVert _{L^2(\R,\mathrm
{d}x)}^2
\\
&&\qquad \leq \max_{x\in\{\phi(0),\phi(0)(1-\alpha\mu)^{-1}\}}\bigl\llvert 1-x\llVert u\rrVert
_{L^2(\R,\mrr x)}^2\bigr\rrvert +\frac{2\phi(0)\alpha\mu}{1-\alpha\mu}\llVert u\rrVert
_{L^2(\R,\mrr x)}^2.\nonumber
\end{eqnarray}
%
The claim follows by \eqref{eq:thm41}, \eqref{eq:secondmom2} and the
rightmost inequality in \eqref{eq:disintens}.\quad\qed
\end{longlist}\noqed
\end{pf*}

\begin{pf*}{Proof of Theorem~\ref{thm:nonlinGaussappr}}
Without loss of generality, we may assume $u\in L^2(\R,\mathrm
{d}x)\cap L^3(\R,\mrr x)\cap L^4(\R,\mrr x)$
(otherwise the claim trivially holds). By the triangular inequality, we have
\begin{eqnarray}
d_W\bigl(\delta_a(u),Z\bigr)\leq d_W
\bigl(\delta_a(u),\delta(u)\bigr)+d_W\bigl(\delta (u),Z
\bigr).
\nonumber
\end{eqnarray}
So, due to Theorem~\ref{thm:nonlinGauss}, we only need to prove
\begin{equation}
\label{eq:dWa} d_W\bigl(\delta_a(u),\delta(u)\bigr)\leq
\frac{2\phi(0)\alpha\mu}{1-\alpha
\mu}\llVert u\rrVert _{L^1(\R,\mrr x)}.
\end{equation}
We have
\begin{eqnarray}
d_W\bigl(\delta_a(u),\delta(u)\bigr)&=&\sup
_{h\in\operatorname{Lip}(1)}\bigl\llvert \mathrm {E}\bigl[h\bigl(\delta_a(u)
\bigr)\bigr]-\mathrm{E}\bigl[h\bigl(\delta(u)\bigr)\bigr]\bigr\rrvert
\nonumber
\\
&\leq&\mathrm{E}\bigl[\bigl\llvert \delta_a(u)-\delta(u)\bigr\rrvert
\bigr] 
=\mathrm{E} \biggl[\biggl\llvert \int
_{\R}u(t)\lambda(t) \,\mathrm {d}t-\widehat{\lambda}\int
_{\R}u(t) \,\mrr t\biggr\rrvert \biggr]
\nonumber
\\
&\leq&\int_{\R}\bigl\llvert u(t)\bigr\rrvert \mathrm{E}\bigl[
\bigl\llvert \lambda(t)-\widehat{\lambda}\bigr\rrvert \bigr] \,\mrr t.
\nonumber
\end{eqnarray}
Inequality \eqref{eq:dWa} then follows by bounding the term
$\mathrm{E}[\llvert  \lambda(t)-\widehat{\lambda}\rrvert  ]$, $t\in\R$, with
the quantity $2\phi(0)\alpha\mu/(1-\alpha\mu)$ (since $\widehat
{\lambda}\in[\phi(0),\phi(0)(1-\alpha\mu)^{-1}]$
the same arguments for \eqref{eq:L1lambda} work).
\end{pf*}

\section{The case of stationary linear Hawkes processes}\label{sec:linear}

Let $N$ be a stationary linear Hawkes process with parameters $(\nu
,h)$ and $\mu:=\int_0^\infty h(t) \,\mrr t<1$. Taking the mean
of its stochastic intensity
we easily see that the intensity of $N$ is equal to
\begin{equation}
\label{eq:intensity} \lambda=\frac{\nu}{1-\mu}
\end{equation}
and so the ``approximated'' first chaos reads
\[
\delta_a(u):=\int_{\R}u(t) \bigl(N(\mrr t)-\nu(1-\mu )^{-1}\,\mrr t\bigr).
\]

\subsection{Explicit Gaussian bounds for the first chaos of linear Hawkes processes}\label{subsec:5/1}

The knowledge of the intensity allows to improve the bounds \eqref
{eq:thm41explicit} and \eqref{eq:thm41explicitappr}
specialized to the linear case. More precisely, the following theorem holds.

\begin{Theorem}\label{thm:linGauss1}
Assume $h:\R_+\to[0,\infty)$ locally bounded and $\mu<1$. Let $N$
be a stationary linear Hawkes process with parameters $(\nu,h)$.
If $u\in L^1(\R,\mrr x)$, then
\begin{eqnarray}
d_W\bigl(\delta(u),Z\bigr)&\leq&\mathfrak{L}\label{eq:thm41explicitlin}
\end{eqnarray}
and
\begin{eqnarray}
d_W\bigl(\delta_a(u),Z\bigr)&\leq&\mathfrak{L}
+
\frac{2\nu\mu}{1-\mu}\llVert u\rrVert _{L^1(\R,\mrr x)},\label
{eq:thm41explicitapprlin}
\end{eqnarray}
where
\begin{eqnarray}
\mathfrak{L}&:=&\sqrt{2/\pi}\biggl\llvert 1-\frac{\nu}{1-\mu}\llVert u\rrVert
_{L^2(\R
,\mrr x)}^2\biggr\rrvert +\frac{\nu}{1-\mu}\llVert u\rrVert
_{L^3(\R,\mathrm
{d}x)}^3
\nonumber
\\
&&{}+\frac{2\sqrt{2/\pi}\nu\mu(2-\mu)}{(1-\mu)^2}\llVert u\rrVert _{L^2(\R
,\mrr x)}^2+
\frac{\nu\mu}{(1-\mu)^2}\llVert u\rrVert _{L^2(\R,\mathrm
{d}x)}\bigl\llVert u^2
\bigr\rrVert _{L^2(\R,\mrr x)}.
\nonumber
\end{eqnarray}
\end{Theorem}

In the linear case, alternative explicit Gaussian bounds may be
obtained using the spectral theory of self-exciting processes;
see \cite{hawkes}. See also \cite{daley}, pages~\mbox{303--309}.

\begin{Theorem}\label{thm:linGaussappvecchio}
Under assumptions and notation of Theorem~\ref{thm:linGauss1}, if
moreover $h\in L^2(\R_+,\mrr x)$, then
\begin{eqnarray}
d_W\bigl(\delta(u),Z\bigr)&\leq&\mathfrak{L}'\label{eq:cov1}
\end{eqnarray}
and
\begin{eqnarray}\label{eq:cov2}
&& d_W\bigl(\delta_a(u),Z\bigr)
\nonumber\\[-8pt]\\[-8pt]\nonumber
&&\qquad \leq \mathfrak{L}'+
\frac{\sqrt{\nu}}{(1-\mu)^{3/2}} \min\bigl\{\mu\llVert u\rrVert _{L^2(\R,\mrr x)},\llVert h\rrVert
_{L^2(\R_+,\mathrm
{d}x)}\llVert u\rrVert _{L^1(\R,\mrr x)}\bigr\},\nonumber
\end{eqnarray}
where
\begin{eqnarray*}
\mathfrak{L}'&:=&\sqrt{2/\pi}
\biggl( \biggl(1-
\frac{\nu}{1-\mu}\llVert u\rrVert _{L^2(\R,\mrr x)}^2
\biggr)^2
\\
&&{} + \frac{\nu}{(1-\mu)^3}\min\bigl\{\mu^2\bigl\llVert
u^2\bigr\rrVert _{L^2(\R,\mathrm
{d}x)}^2,\llVert h\rrVert
_{L^2(\R_+,\mrr x)}^2\bigl\llVert u^2\bigr\rrVert
_{L^1(\R,\mathrm
{d}x)}^2\bigr\}\biggr)^{1/2}
\nonumber
\\
&&{}+\frac{\nu}{1-\mu}\llVert u\rrVert _{L^3(\R,\mrr x)}^3+
\frac{2\sqrt
{2/\pi}\nu\mu}{(1-\mu)^2}\llVert u\rrVert _{L^2(\R,\mrr x)}^2
\\
&&{} +
\frac{\nu\mu}{(1-\mu)^2}\llVert u\rrVert _{L^2(\R,\mrr x)}\bigl\llVert u^2
\bigr\rrVert _{L^2(\R,\mrr x)}.
\nonumber
\end{eqnarray*}
\end{Theorem}

The proofs of Theorems~\ref{thm:linGauss1} and~\ref
{thm:linGaussappvecchio} are given in Section~\ref{sec:proofs2}.

Next proposition (whose proof is a simple consequence of the elementary
inequality $\sqrt{a^2+b^2}\leq\llvert  a\rrvert  +\llvert  b\rrvert  $, $a,b\in\R$, and therefore omitted)
provides sufficient conditions under which the bounds of Theorem~\ref
{thm:linGaussappvecchio}
improve the bounds of Theorem~\ref{thm:linGauss1}. Hereafter, for ease
of notation,
we denote by $\widetilde{\mathfrak{L}}$ the right-hand side of~\eqref
{eq:thm41explicitapprlin} and by $\widetilde{\mathfrak{L}'}$ the
right-hand side of
\eqref{eq:cov2}.

\begin{Proposition}\label{prop:compare}
Under assumptions and notation of Theorem~\ref
{thm:linGaussappvecchio}, we have:

\begin{longlist}[(ii)]
\item[(i)] If
\begin{equation}
\label{eq:compare} \nu\geq\frac{1}{4(1-\mu)}\min \biggl\{\frac{\llVert  u^2\rrVert  _{L^2(\R
,\mrr x)}^2}{\llVert  u\rrVert  _{L^2(\R,\mrr x)}^4},
\frac{\llVert  h\rrVert  _{L^2(\R_+,\mrr x)}^2}{\mu^2} \biggr\} 
\end{equation}
then
$\mathfrak{L}'\leq\mathfrak{L}$.

\item[(ii)] If
%
\begin{eqnarray}\label{eq:compare1}
\nu &\geq& \frac{1}{4(1-\mu)}\max \biggl\{ \min \biggl\{
\frac{\llVert  u^2\rrVert  _{L^2(\R,\mrr x)}^2}{\llVert  u\rrVert  _{L^2(\R
,\mrr x)}^4}, \frac{\llVert  h\rrVert  _{L^2(\R_+,\mrr x)}^2}{\mu^2} \biggr\},
\nonumber\\[-8pt]\\[-8pt]\nonumber
&&{} \min \biggl\{\frac{\llVert  u\rrVert  _{L^2(\R,\mrr x)}^2}{\llVert  u\rrVert  _{L^1(\R
,\mrr x)}^2},\frac{\llVert  h\rrVert  _{L^2(\R_+,\mrr x)}^2}{\mu
^2} \biggr\} \biggr\} 
\end{eqnarray}
then $\widetilde{\mathfrak{L}'}\leq\widetilde{\mathfrak{L}}$.
\end{longlist}
\end{Proposition}

\subsection{A quantitative central limit theorem for linear Hawkes
processes}\label{subsec:5/2}

The following quantitative central limit theorem in the Wasserstein
distance is an immediate consequence of Theorems
\ref{thm:linGauss1} and~\ref{thm:linGaussappvecchio}.

\begin{Corollary}\label{cor:QCLT}
For $\varepsilon>0$, assume $h_\varepsilon:\R_+\to[0,\infty)$
locally bounded functions and such that
$\mu_\varepsilon:=\int_0^\infty h_\varepsilon(x) \,\mrr x<1$.
Let $N_\varepsilon$ be a stationary linear Hawkes process with parameters
$(\nu_\varepsilon,h_\varepsilon)$ and take $u_\varepsilon\in L^1(\R
,\mrr x)$. Then:

\begin{longlist}[(ii)]
\item[(i)]
\begin{eqnarray} \label{eq:bound1}
d_W\bigl(\delta^{(\varepsilon)}(u_\varepsilon),Z\bigr)&\leq&\min
\bigl\{\mathfrak {L}_\varepsilon,\mathfrak{L}_{\varepsilon}'
\ind_{\{h_\varepsilon
\in L^2(\R_+,\mrr x)\}}\bigr\} +\mathfrak{L}_\varepsilon\ind_{\{h_\varepsilon\notin L^2(\R
_+,\mrr x)\}},
\nonumber\\[-8pt]\\[-8pt]
\eqntext{\varepsilon>0}
\end{eqnarray}
and
\begin{eqnarray}\label{eq:bound2}
d_W\bigl(\delta_a^{(\varepsilon)}(u_\varepsilon),Z
\bigr)&\leq&\min\bigl\{ \widetilde{\mathfrak{L}_{\varepsilon}},\widetilde{
\mathfrak {L}_{\varepsilon}'} \ind_{\{h_\varepsilon\in L^2(\R_+,\mrr x)\}}\bigr\}+
\widetilde{\mathfrak{L}_{\varepsilon}}\ind_{\{h_\varepsilon\notin
L^2(\R_+,\mrr x)\}},
\nonumber\\[-8pt]\\[-8pt]
\eqntext{\varepsilon>0.}
\end{eqnarray}
Here, $\mathfrak{L}_\varepsilon$, $\mathfrak{L}_{\varepsilon}'$,
$\widetilde{\mathfrak{L}_{\varepsilon}}$ and $\widetilde{\mathfrak
{L}_{\varepsilon}'}$
are defined as $\mathfrak{L}$, $\mathfrak{L}'$, $\widetilde
{\mathfrak{L}}$ and $\widetilde{\mathfrak{L}'}$, respectively, with
$\nu_\varepsilon$, $\mu_\varepsilon$,
$u_\varepsilon$ and $h_\varepsilon$ in place of $\nu$, $\mu$, $u$
and $h$, respectively;
\begin{eqnarray*}
\delta^{(\varepsilon)}(u_\varepsilon)&:=&\int_{\R}u_\varepsilon
(t) \bigl(N_\varepsilon(\mrr t)-\lambda_\varepsilon(t) \,\mrr t
\bigr),
\\
\lambda_\varepsilon(t) &:=& \nu_\varepsilon+\int
_{(-\infty
,t)}h_\varepsilon(t-s)N_\varepsilon(\mrr s),
\\
\delta_a^{(\varepsilon)}(u_\varepsilon)&:=&\int
_{\R}u_\varepsilon (t) \bigl(N_\varepsilon(\mrr t)-\nu_\varepsilon(1-\mu_\varepsilon )^{-1}
\,\mrr t\bigr).
\end{eqnarray*}

\item[(ii)] If, as $\varepsilon\to0$,
\begin{eqnarray}
\label{eq:lim1} \mu_\varepsilon&\to&0,
\\
\label{eq:lim2} \nu_\varepsilon\llVert u_\varepsilon\rrVert
_{L^2(\R,\mrr x)}^2&\to&1,
\\
\label{eq:lim3} \nu_\varepsilon\llVert u_\varepsilon\rrVert
_{L^3(\R,\mrr x)}^3&\to&0,
\\
\label{eq:lim4} \nu_\varepsilon(\mu_\varepsilon)^2\bigl
\llVert (u_\varepsilon)^2\bigr\rrVert _{L^2(\R
,\mrr x)}^2 &\to& 0,
\end{eqnarray}
then
\[
d_W\bigl(\delta^{(\varepsilon)}(u_\varepsilon),Z
\bigr)\to0,\qquad\mbox{as }\varepsilon\to0.
\]
If moreover, as $\varepsilon\to0$,
\begin{equation}
\label{eq:lim4ter} \nu_\varepsilon\mu_\varepsilon\llVert u_\varepsilon
\rrVert _{L^1(\R,\mathrm
{d}x)}\to0, 
\end{equation}
then
\[
d_W\bigl(\delta_a^{(\varepsilon)}(u_\varepsilon),Z
\bigr)\to0,\qquad\mbox{as }\varepsilon\to0.
\]
This latter limit holds even if we replace condition \eqref
{eq:lim4ter} with
\begin{equation}
\label{eq:SPMS1} h_\varepsilon\in L^2(\R_+,\mrr x),\qquad
\varepsilon>0. 
\end{equation}
%
\end{longlist}
\end{Corollary}

%
\begin{Remark}
In this remark, we compare Corollary~\ref{cor:QCLTnonlinear},
specialized to the case
of a self-exciting process $N_\varepsilon$ with parameters $(\nu
_\varepsilon,h_\varepsilon)$, with Corollary
\ref{cor:QCLT}. First, we note that the upper bounds \eqref
{eq:bound1} and \eqref{eq:bound2} improve the upper bounds
\eqref{eq:bound1bis} and \eqref{eq:bound2bis}, respectively. Second,
we note that conditions \eqref{eq:lim1}--\eqref{eq:lim4ter}
coincide with conditions \eqref{eq:lim1bis}--\eqref{eq:lim5bis}.
Finally, we note that
in Corollary~\ref{cor:QCLT} we deduce the convergence to zero
of the family $\{d_W(\delta_a^{(\varepsilon)}(u_\varepsilon),Z)\}
_{\varepsilon>0}$ even replacing condition \eqref{eq:lim4ter} with
the alternative condition \eqref{eq:SPMS1}.
\end{Remark}

We conclude this subsection with an example.

\begin{Example}\label{ex:2}
Let $I_\varepsilon$, $\varepsilon>0$, be a given family of bounded
Borel sets, $I_\varepsilon$ with Lebesgue measure $\ell_\varepsilon$,
and $\nu_\varepsilon>0$, $\varepsilon>0$, be a family of positive
constants. Let $\mu_\varepsilon$, $\varepsilon>0$,
be a collection of positive numbers such that $\mu_\varepsilon<1$,
$\varepsilon>0$, and define the functions
$h_\varepsilon(t):=\mu_\varepsilon f_\varepsilon(t)$, $\varepsilon
>0$, $t>0$, where
$f_\varepsilon$ is a locally bounded probability density $($with
respect to the Lebesgue measure$)$
on $(0,\infty)$ such that $f_\varepsilon\in L^2(\R_+,\mrr x)$,
$\varepsilon>0$. Hereafter, we consider
the family $N_\varepsilon$, $\varepsilon>0$, of stationary linear
Hawkes processes with parameters
$(\nu_\varepsilon,h_\varepsilon)$, $\varepsilon>0$,
and the functions
\[
u_\varepsilon(t):=\frac{1}{\sqrt{\vafrac{\nu_\varepsilon\ell
_\varepsilon}{1-\mu_\varepsilon}}}\ind_{I_\varepsilon}(t),\qquad
\varepsilon>0, t\in\R.
\]
Using the expressions of the $L^p$-norms of $u_\varepsilon$ computed
in the Example~\ref{ex:1} [clearly setting $\alpha_\varepsilon=1$
and $\phi_\varepsilon(0)=\nu_\varepsilon$ therein],
one may easily see that conditions \eqref{eq:compare} and \eqref
{eq:compare1} are both equivalent to
\begin{equation}
\label{eq:ineqbestbound} \nu_\varepsilon\geq\frac{1}{4(1-\mu_\varepsilon)}\min\bigl\{\ell
_\varepsilon^{-1},\llVert f_\varepsilon\rrVert
_{L^2(\R_+,\mrr x)}^2\bigr\},\qquad\varepsilon>0.
\end{equation}
Note\vspace*{1pt} also that the square integrability of $f_\varepsilon$ implies
$h_\varepsilon\in L^2(\R_+,\mrr x)$. Therefore, under \eqref
{eq:ineqbestbound}, by Proposition~\ref{prop:compare} we deduce
$\mathfrak{L}_\varepsilon'\leq\mathfrak{L}_\varepsilon$
and $\widetilde{\mathfrak{L}_\varepsilon'}\leq\widetilde{\mathfrak
{L}_\varepsilon}$, $\varepsilon>0$.
So, under
\eqref{eq:ineqbestbound}, by Corollary~\ref{cor:QCLT} we have
%
\begin{eqnarray*}
&& d_W\bigl(\delta^{(\varepsilon)}(u_\varepsilon),Z\bigr)
\\
&&\qquad \leq \mathfrak {L}_\varepsilon' = \bigl(\sqrt{2/\pi}\min\bigl\{(
\ell_\varepsilon)^{-1/2},\llVert f_\varepsilon\rrVert
_{L^2(\R_+,\mrr x)}\bigr\}+(\ell_\varepsilon)^{-1/2}\bigr)
\frac{\mu_\varepsilon}{\sqrt{\nu_\varepsilon(1-\mu_\varepsilon)}}
\\
&&\quad\qquad{}+\sqrt{\frac{1-\mu_\varepsilon}{\nu_\varepsilon\ell_\varepsilon
}}+2
\sqrt{2/\pi}\frac{\mu_\varepsilon}{1-\mu_\varepsilon},\qquad \varepsilon>0
\nonumber
\end{eqnarray*}
and
\begin{eqnarray*}
&& d_W\bigl(\delta_a^{(\varepsilon)}(u_\varepsilon),Z
\bigr)
\\
&&\qquad \leq\widetilde {\mathfrak{L}_\varepsilon'}
=  \bigl(\sqrt{2/
\pi}\min\bigl\{(\ell_\varepsilon)^{-1/2},\llVert f_\varepsilon
\rrVert _{L^2(\R_+,\mrr x)}\bigr\}+(\ell_\varepsilon)^{-1/2}\bigr)
\frac{\mu_\varepsilon}{\sqrt{\nu_\varepsilon(1-\mu_\varepsilon
)}}
\\
&&\quad\qquad {}+\sqrt{\frac{1-\mu_\varepsilon}{\nu_\varepsilon\ell_\varepsilon
}}+\bigl(2
\sqrt{2/\pi}+\min\bigl\{1,\sqrt{\ell_\varepsilon}\llVert f_\varepsilon
\rrVert _{L^2(\R_+,\mrr x)}\bigr\}\bigr) \frac{\mu_\varepsilon}{1-\mu_\varepsilon},\qquad\varepsilon >0.
\nonumber
\end{eqnarray*}
Finally, one easily sees that if
\begin{equation}
\label{eq:ultimahyp} \mu_\varepsilon\to0\quad\mbox{and}\quad\nu_\varepsilon
\ell _\varepsilon\to\infty,\qquad\mbox{as }\varepsilon\to0
\end{equation}
then
\[
d_W\bigl(\delta^{(\varepsilon)}(u_\varepsilon),Z\bigr)\to0\quad
\mbox {and}\quad d_W\bigl(\delta_a^{(\varepsilon)}(u_\varepsilon),Z
\bigr)\to 0,\qquad\mbox{as }\varepsilon\to0.
\]
It has to be noticed that a straightforward computation shows that
condition \eqref{eq:ultimahyp} implies conditions \eqref
{eq:lim1}--\eqref{eq:lim4}, but
does not imply condition \eqref{eq:lim4ter}.
\end{Example}

\subsection{Proofs of Theorems \texorpdfstring{\protect\ref{thm:linGauss1}}{5.1} and
\texorpdfstring{\protect\ref{thm:linGaussappvecchio}}{5.2}}\label{sec:proofs2}

\mbox{}

\begin{pf*}{Proof of Theorem~\ref{thm:linGauss1}} The
claim follows by an obvious modification of the proofs of Theorems~\ref
{thm:nonlinGauss} and~\ref{thm:nonlinGaussappr}.
For instance, to get \eqref{eq:thm41explicitlin} it suffices to modify
the proof of Theorem~\ref{thm:nonlinGauss} as follows. We combine
inequality \eqref{eq:thm41} [taking~therein $\lambda=\nu(1-\mu
)^{-1}$ and $\alpha=1$]
with inequalities \eqref{eq:secondmom} and \eqref{eq:L1lambda}
[taking therein $\lambda=\nu(1-\mu)^{-1}$, $\alpha=1$ and $\phi
(0)=\nu$].
Note that, due to the knowledge of~$\lambda$, we do not need anymore
to further bound
the quantity $\llvert  1-\lambda\llVert  u\rrVert  _{L^2(\R,\mrr x)}^2\rrvert  $ as in~\eqref
{eq:secondmom2}.
\end{pf*}

\begin{pf*}{Proof of Theorem~\ref{thm:linGaussappvecchio}}
The claim is clearly true if $u\notin L^2(\R
,\mrr x)\cap L^3(\R,\mrr x)\cap L^4(\R,\mrr x)$.
So we shall assume these integrability conditions. We first prove the
bound \eqref{eq:cov1}.
By \eqref{eq:thm41} [with $\lambda=\nu(1-\mu)^{-1}$ and $\alpha
=1$], the Cauchy--Schwarz inequality and the stationarity of $N$,
we have
\begin{eqnarray}\label{eq:thm42}
d_W\bigl(\delta(u),Z\bigr)&\leq&\sqrt{2/\pi}\sqrt{1-2
\frac{\nu}{1-\mu}\llVert u\rrVert _{L^2(\R,\mrr x)}^2+\int
_{\R^2}\bigl\llvert u(t)u(s)\bigr\rrvert ^2
\mathrm{E}\bigl[\lambda (t)\lambda(s)\bigr] \,\mrr t\,\mrr s}
\nonumber
\\
&&{}+\frac{\nu}{1-\mu}\llVert u\rrVert _{L^3(\R,\mrr x)}^3
+\frac{2\nu\mu}{(1-\mu)^2}\sqrt{2/\pi}\llVert u\rrVert _{L^2(\R,\mrr x)}^{2}
\\
&&{} +\frac{\nu\mu}{(1-\mu)^2}\llVert u\rrVert _{L^2(\R,\mrr x)}\bigl\llVert u^2
\bigr\rrVert _{L^2(\R,\mrr x)}.\nonumber
\end{eqnarray}
%
By \eqref{eq:intensity} and again the stationarity of $N$, we deduce
\begin{eqnarray}
\mathrm{E}\bigl[\lambda(t)\lambda(s)\bigr] &=& \lambda^2+\operatorname{Cov}
\bigl(\lambda (t),\lambda(s)\bigr) 
\nonumber\\[-8pt]\\[-8pt]\nonumber
&=& \biggl(\frac{\nu}{1-\mu} \biggr)^2+\operatorname{Cov} \biggl(\int
_{\R }h_t(u)N(\mrr u),\int
_{\R}h_s(u)N(\mrr u) \biggr),\label{eq:tromb1}
\end{eqnarray}
where\vspace*{1pt} we set $h_t(u):=\ind_{(-\infty,t)}(u)h(t-u)$. In the following,
for $f\in L^1(\R,\mrr x)\cap L^2(\R,\mrr x)$,
we denote by $\widehat{f}(\omega):=\int_{\R}\mathrm{e}^{i\omega
t}f(t) \,\mrr t$
the Fourier transform of $f$, and we extend the definition of $h$ on
$(-\infty,0]$ setting $h(t):=0$ for $t\leq0$.
By the results in \cite{hawkes}, we have (see also formulas $(8)$ and
$(24)$ in \cite{bremaudmassoulie})
\begin{eqnarray}\label{eq:tromb2}
&& \operatorname{Cov} \biggl(\int_{\R}h_t(u)N(\mrr u),\int_{\R
}h_s(u)N(\mrr u)
\biggr)
\nonumber\\[-8pt]\\[-8pt]\nonumber
&&\qquad = \frac{\nu}{2\pi(1-\mu)}\int_{\R}\widehat{h_t}(
\omega)\widehat {h_s}(\omega)\frac{1}{\llvert  1-\widehat{h}(\omega)\rrvert  ^2} \,\mrr \omega.\nonumber
\end{eqnarray}
Note that
\begin{equation}
\label{eq:tromb3} \bigl\llvert 1-\widehat{h}(\omega)\bigr\rrvert \geq\big\llvert 1-
\big\rrvert \widehat{h}(\omega)\big\llvert \big\rrvert \geq 1-\bigl\llvert \widehat{h}(\omega)
\bigr\rrvert \geq1-\mu>0,\qquad\omega\in\R
\end{equation}
and that $\widehat{h_t}(\omega)=\mathrm{e}^{i\omega t}\widehat
{h}(-\omega)$ (since $h$ has a positive support).
Therefore,
\begin{eqnarray}
&&\int_{\R^2}\bigl\llvert u(t)u(s)\bigr\rrvert ^2
\mathrm{E}\bigl[\lambda(t)\lambda(s)\bigr] \,\mathrm {d}t\,\mrr s
\nonumber
\\
&&\qquad = \biggl(\frac{\nu}{1-\mu} \biggr)^2\llVert u\rrVert
_{L^2(\R,\mrr x)}^4
\nonumber\\[-8pt]\label{eq:FubiniFourier} \\[-8pt]\nonumber
&&\quad\qquad{} + \frac{\nu}{2\pi(1-\mu)}\int_{\R}
\biggl(\int_{\R
}\bigl\llvert u(t)\bigr\rrvert ^2
\widehat{h_t}(\omega) \,\mrr t \biggr)^2
\frac
{1}{\llvert  1-\widehat{h}(\omega)\rrvert  ^2} \,\mrr \omega
\\
&&\qquad  \leq \biggl(\frac{\nu}{1-\mu} \biggr)^2\llVert u\rrVert
_{L^2(\R,\mrr x)}^4+ \frac{\nu}{2\pi(1-\mu)^3}\int_{\R}
\biggl\llvert \int_{\R
}\bigl\llvert u(t)\bigr\rrvert
^2\widehat{h_t}(\omega) \,\mrr t\biggr\rrvert
^2 \,\mathrm {d}\omega.\label{eq:FubiniFourier1}
\end{eqnarray}
%
In \eqref{eq:FubiniFourier}, we used Fubini's theorem, which is
applicable since
\begin{eqnarray*}
&& \int_{\R^3}\bigl\llvert u(t)u(s)\bigr\rrvert ^2
\bigl\llvert \widehat{h_t}(\omega)\widehat{h_s}(\omega )
\bigr\rrvert \frac{1}{\llvert  1-\widehat{h}(\omega)\rrvert  ^2} \,\mrr s\,\mathrm {d}t\,\mrr \omega
\\
&&\qquad \leq
\frac{1}{(1-\mu)^2}\llVert u\rrVert _{L^2(\R,\mrr x)}^4\int
_{\R
}\bigl\llvert \widehat{h}(-\omega)\bigr\rrvert
^2 \,\mrr \omega
\\
&&\qquad =\frac{2\pi}{(1-\mu)^2}\llVert u\rrVert _{L^2(\R,\mrr x)}^4\llVert h
\rrVert _{L^2(\R
,\mrr x)}^2<\infty,
\end{eqnarray*}
where in the latter equality we used Parseval's identity. Setting
$\check{f}(x):=f(-x)$, $u^2(\cdot):=u(\cdot)^2$ and letting the
symbol $*$ denote the convolution product, we have
\[
\int_{\R}\bigl\llvert u(t)\bigr\rrvert ^2
\widehat{h_t}(\omega) \,\mrr t=\widehat {h}(-\omega)
\widehat{u^2}(\omega)=\widehat{\check{h}}(\omega )
\widehat{u^2}(\omega) =\widehat{\check{h}*u^2}(\omega).
\]
Consequently, using again the Parseval identity, we deduce
\begin{eqnarray}
\int_{\R}\biggl\llvert \int_{\R}\bigl
\llvert u(t)\bigr\rrvert ^2\widehat{h_t}(\omega) \,\mathrm {d}t\biggr\rrvert ^2 \,\mrr \omega&=&\int_{\R}
\bigl\llvert \widehat{\check {h}*u^2}(\omega)\bigr\rrvert
^2 \,\mrr \omega 
=2\pi\bigl\llVert
\check{h}*u^2\bigr\rrVert _{L^2(\R,\mrr x)}^2.
\nonumber
\end{eqnarray}
By this relation, \eqref{eq:FubiniFourier1} and standard properties of
the convolution (see, e.g., Theorem IV.15 in \cite{brezis}), we have
\begin{eqnarray*}
&& \int_{\R^2}\bigl\llvert u(t)u(s)\bigr\rrvert ^2
\mathrm{E}\bigl[\lambda(t)\lambda(s)\bigr] \,\mathrm {d}t\,\mrr s
\\
&&\qquad \leq
\frac{\nu^2}{(1-\mu)^2}\llVert u\rrVert _{L^2(\R,\mrr x)}^4+
\frac{\nu}{(1-\mu)^3}\bigl\llVert \check{h}*u^2\bigr\rrVert
_{L^2(\R,\mathrm
{d}x)}^2
\nonumber
\\
&&\qquad \leq\frac{\nu^2}{(1-\mu)^2}\llVert u\rrVert _{L^2(\R,\mrr x)}^4
\\
&&\quad\qquad{}+
\frac{\nu}{(1-\mu)^3}\min\bigl\{\mu^2\bigl\llVert u^2\bigr
\rrVert _{L^2(\R,\mathrm
{d}x)}^2,\llVert h\rrVert _{L^2(\R,\mrr x)}^2
\bigl\llVert u^2\bigr\rrVert _{L^1(\R,\mathrm
{d}x)}^2\bigr\}.
\nonumber
\end{eqnarray*}
The claim follows combining this inequality with \eqref{eq:thm42}. We
now prove the bound \eqref{eq:cov2}.
By the triangular inequality and \eqref{eq:cov1}, we only need to prove
\[
d_W\bigl(\delta_a(u),\delta(u)\bigr)\leq
\frac{\sqrt{\nu}}{(1-\mu
)^{3/2}}\min\bigl\{\mu\llVert u\rrVert _{L^2(\R,\mrr x)},\llVert h\rrVert
_{L^2(\R,\mathrm
{d}x)}\llVert u\rrVert _{L^1(\R,\mrr x)}\bigr\}.
\]
Note that
\begin{eqnarray}
d_W\bigl(\delta_a(u),\delta(u)\bigr)&=&\sup
_{h\in\operatorname{Lip}(1)}\bigl\llvert \mathrm {E}\bigl[h\bigl(\delta_a(u)
\bigr)\bigr]-\mathrm{E}\bigl[h\bigl(\delta(u)\bigr)\bigr]\bigr\rrvert
\nonumber
\\
&\leq&\mathrm{E}\bigl[\bigl\llvert \delta_a(u)-\delta(u)\bigr\rrvert
\bigr] 
=\mathrm{E} \biggl[\biggl\llvert \int
_{\R}u(t)\lambda(t) \,\mathrm {d}t-\lambda\int
_{\R}u(t) \,\mrr t\biggr\rrvert \biggr]
\nonumber
\end{eqnarray}
and so by the Cauchy--Schwarz inequality we have
\begin{eqnarray} \label{eq:FubiniSinga}
&& d_W\bigl(\delta_a(u),\delta(u)\bigr)\nonumber
\\
&&\qquad \leq \biggl(
\mathrm{E} \biggl[\biggl\llvert \int_{\R}u(t)\lambda(t)
\,\mrr t-\frac{\nu}{1-\mu}\int_{\R}u(t) \,\mrr t\biggr
\rrvert ^2 \biggr] \biggr)^{1/2}
\\
&&\qquad = \sqrt{\int_{\R^2}u(t)u(s)\mathrm{E}\bigl[\lambda(t)
\lambda(s)\bigr] \,\mrr t\,\mrr s- \biggl(\frac{\nu}{1-\mu}
\biggr)^2 \biggl(\int_{\R}u(t) \,\mrr t
\biggr)^2}.\nonumber
\end{eqnarray}
%
To upper bound the first addend inside the square root,
we repeat the arguments above. So, by \eqref{eq:tromb1}, \eqref
{eq:tromb2}, \eqref{eq:tromb3},
Fubini's theorem and Parseval's identity,
we have
\begin{eqnarray}\label{eq:FubiniSinga2}
&& \int_{\R^2}u(s)u(t)\mathrm{E}\bigl[\lambda(s)\lambda(t)
\bigr] \,\mathrm {d}s\,\mrr t\nonumber
\\
&&\qquad \leq \biggl(\frac{\nu}{1-\mu} \biggr)^2
\biggl(\int_{\R}u(t) \,\mrr t \biggr)^2+
\frac{\nu}{2\pi(1-\mu)^3}\int_{\R}\biggl\llvert \int
_{\R}u(t)\widehat {h_t}(\omega) \,\mrr t\biggr
\rrvert ^2 \,\mrr \omega\nonumber
\\
&&\qquad = \biggl(\frac{\nu}{1-\mu} \biggr)^2 \biggl(\int
_{\R}u(t) \,\mathrm {d}t \biggr)^2+
\frac{\nu}{(1-\mu)^3}\llVert \check{h}*u\rrVert _{L^2(\R,\mathrm
{d}x)}^2
\\
&&\qquad \leq\biggl(\frac{\nu}{1-\mu} \biggr)^2 \biggl(\int
_{\R}u(t) \,\mrr t \biggr)^2\nonumber
\\
&&\qquad\quad{} +
\frac{\nu}{(1-\mu)^3}\min\bigl\{\mu^2\llVert u\rrVert
_{L^2(\R,\mrr x)}^2,\llVert h\rrVert _{L^2(\R,\mrr x)}^2
\llVert u\rrVert _{L^1(\R,\mrr x)}^2\bigr\}.\nonumber
\end{eqnarray}
Note that in \eqref{eq:FubiniSinga} we used Fubini's theorem to
exchange the double integral with the expectation. This is justified by
the fact that
inequality \eqref{eq:FubiniSinga2} holds replacing $u$ with $\llvert  u\rrvert  $
and the resulting right-hand side is finite. The proof is complete.
\end{pf*}

\section*{Acknowledgments}
The author thanks the Editor and two anonymous reviewers for a careful
reading of the paper.



\printaddresses
\end{document}